\documentclass[11pt,a4]{article}
\usepackage{amsmath,amssymb,thc,amsfonts,dsfont,enumerate,a4wide}
\usepackage{hyperref}
\usepackage{color}

\theoremstyle{change}
\theorembodyfont{\itshape}
\newtheorem{theorem}{\bf Theorem}[section]
\newtheorem{proposition}[theorem]{\bf Proposition}
\newtheorem{lemma}[theorem]{\bf Lemma}
\newtheorem{corollary}[theorem]{\bf Corollary}

\theorembodyfont{\rm}

\newtheorem{remark}[theorem]{\bf Remark}
\newtheorem{remarks}[theorem]{\bf Comments}
\newtheorem{example}[theorem]{\bf Example}

\newcommand{\conver}{\mathop{\longrightarrow}}
\newcommand{\inprob}{\ \conver_{\cP}\ }
\newcommand{\indist}{\longrightarrow_{\mathcal{D}}}

\newcommand{\cN}{{\mathcal N}}
\newcommand{\cP}{{\mathcal P}}

\newcommand{\cB}{{\mathcal B}}
\newcommand{\cF}{{\mathcal F}}

\newcommand{\cS}{{\mathcal S}}

\newcommand{\GN}{{\mathds{N}}}
\newcommand{\GR}{{\mathds{R}}}
\newcommand{\GC}{{\mathds{C}}}

\newcommand{\GS}{{\mathds{S}}}

\newcommand{\bE}{{\mathbb{E}}}
\newcommand{\bP}{{\mathbb{P}}}

\newcommand{\bV}{\mathbb{V}}

\newcommand{\Idn}{\mbox{ $1\mskip-1.6\thinmuskip$\textrm{I}}} 
\newcommand{\Idnem}{\mbox{ $1\mskip-1.6\thinmuskip$\text{\em I}}} 
\begin{document}

\title{Stable limits for Markov chains\\ via the Principle of Conditioning}

\author{Mohamed El Machkouri${}^{1}$\footnote{E-mail: mohamed.elmachkouri@univ-rouen.fr}, Adam Jakubowski${}^{2}$\footnote{E-mail: adjakubo@mat.umk.pl} \\ and Dalibor Voln\'y${}^1$\footnote{E-mail: dalibor.volny@univ-rouen.fr} \\[3mm]
${}^1$ Universit\'e de
Rouen Normandie, France\\
${}^{2}$ Nicolaus Copernicus University, Poland}

\date{}

\maketitle

\begin{abstract}
We study limit theorems for partial sums of instantaneous functions of a homogeneous Markov chain on a general state space. The summands are heavy-tailed and the limits are stable distributions. The conditions imposed on the transition operator $P$ of the Markov chain ensure that the limit is the same as if the summands were independent. Such a~scheme admits a physical interpretation, as given in Jara et al. ({\em Ann. Appl. Probab.}, 19 (2009), 2270--2300). 

We considerably extend the results of Jara et al., {\em ibid.} and Cattiaux and Manou-Abi ({\em ESAIM Probab. Stat.}, 18 (2014), 468--486). We show that the theory holds under the assumption of operator uniform integrability in $L^2$ of $P$ (a notion introduced by Wu ({\em J. Funct. Anal.}, 172 (2000), 301--376)) {\em plus} the $L^2$-spectral gap property. If we strengthen the uniform integrability in $L^2$ to the hyperboundedness, then the $L^2$-spectral gap property can be relaxed to the strong mixing at geometric rate (in practice: to geometric ergodicity).

We provide an example of a Markov chain on a countable space that is uniformly integrable in $L^2$ (and admits an $L^2$-spectral gap), while it is not hyperbounded. Moreover, we show by example that hyperboundedness is still a  weaker property than $\phi$-mixing, what enlarges the range of models of interest.

What makes our assumptions working is a new, efficient version of the Principle of Conditioning that operates with conditional characteristic functions rather than predictable characteristics.
\end{abstract}


\noindent{\em MSClassification 2010:} 60F05, 60F17, 60E07, 60J05, 60J35.

\section{Introduction}
Our motivation comes from  the paper by Jara, Komorowski and Olla
 \cite{JKO09}, where  a fractional diffusion was obtained as a scaled 
limit of functionals  of Markov chains forming a probabilistic solution to 
a linear Boltzmann equation. The main tool used in \cite{JKO09}  was 
a functional limit theorem on convergence to stable L\'evy processes 
due to Durret and Resnick \cite{DuRe78} and the assumptions that made 
this functional limit theorem working were $L^2$-spectral gap  and strong
 contractivity properties of the Markov transition operator. In the particular
 example considered in \cite{JKO09} the {\em ultraboundedness} of the
 transition operator was used, but in the general considerations (Theorem 2.4, {\em ibid.}) properties related to a weaker notion of {\em hyperboundedness} 
were assumed (We refer to Section 2 below for formal definitions and 
discussion of all these notions). 

Later Cattiaux and Manou-Abi \cite{CaM-A14} reexamined the limit theorems
 from \cite{JKO09} in the context of the general theory of convergence 
to stable laws for sums of stationary sequences. They considered standard
 mixing conditions ($\phi$-, $\rho$-, $\alpha$- mixing) and anti-clustering
 condition $D'$, introduced in \cite{Dav83} and discussed in \cite{DeJa89} 
(see also \cite{Kriz10}). While the discussion in  \cite{CaM-A14} was quite
 extensive, it did not address the question whether the strong assumption 
of hyperboundedness of the transition operator can be  essentially  weakened.

In the present paper we suggest replacing the hyperboundedness with 
the {\em uniform integrability in $L^2$} ($2$-U.I. in short) of the transition 
operator, a notion  introduced in \cite{Wu00}. We believe that this 
is the proper minimal form 
for operator contractivity whenever limit theorems for Markov chains with
 stable limits are considered. Our main results are formulated in Section \ref{SecResults}. In Theorem \ref{TheMain} we obtain limit theorems assuming the $2$-U.I. condition and the $L^2$-spectral gap property. In Theorem \ref{TheHyper} we assume the hyperboundedness in place of the $2$-U.I. condition, but we weaken the $L^2$-spectral gap property to the geometric ergodicity.
 The proofs of both main results  are
 deferred to Section \ref{SecProofs}. In Section \ref{SecPreliminaries} we gather all necessary information, notation and comments related to the models considered in the paper.

What allows considerable weakening of the assumptions is a new efficient version of 
the Principle of Conditioning that operates with conditional characteristic
 functions rather than predictable characteristics and therefore keeps
integrability requirements at the minimal possible level. Recall that 
the Principle of Conditioning is a heuristic rule that transforms limit theorems 
for independent random variables into limit theorems for  dependent 
random variables. The mentioned above functional limit theorem 
by Durret and Resnick \cite{DuRe78} is a  particular
 manifestation  of this rule.  We state our new result (Theorem \ref{PoCnew}) 
and give  more comments and references on the Principle of Conditioning 
in the Appendix.

In Section \ref{SecExamples} we give four examples, each of different nature.
First we provide an example of a Markov chain with the transition operator 
that is uniformly integrable in $L^2$ (and admits an $L^2$-spectral gap) 
while it is not hypercontractive.  This shows that our theory substantially
 extends that of \cite{JKO09} and   \cite{CaM-A14}.

Then we show that the standard stationary AR(1) sequence with Gaussian
 innovations satisfies the hyperboundedness property and admits an
 $L^2$-spectral gap. It follows that instantaneous functions of this sequence give stable limit theorems without any need of centering in the whole range  $\alpha \in (0,2)$
(and not only for $\alpha \in (0,1)$). This partially answers a conjecture formulated in \cite{Dav83}. 
On the other hand it is well-known that this sequence is not 
$\phi$-mixing what proves that the hyperboundedness is not 
as demanding as it looks like.

Finally we study the problem of $m$-skeletons.  It is known 
that the contraction properties may improve after composition 
of operators. Suppose that some power $P^m$ of the 
transition operator has the desired (by us) properties while
$P^k, \ k=1, 2, \ldots, m-1$ not. It follows that our stable 
limit theorem holds if we sum random variables along
 $m$-skeleton only and the question is whether the limit 
 theorem can be extended to the whole sequence. The 
 answer is ``no" as simple probabilistic examples built upon 
 i.i.d. sequences show. We provide another example, with
  $m=3$, that is more oriented towards thinking in terms 
  of operators.

\section{Preliminaries}\label{SecPreliminaries}

\subsection{Transition operator}
Let $\{X_n\}_{n\geq 0}$ be a Markov chain with state space $(\GS,\cS)$ 
and the transition probability $P(x, dy)$ on $\GS \times \cS$. We will always
 assume that $P(x, dy)$ admits a stationary distribution $\pi$ on 
$(\GS,\cS)$, i.e. 
\begin{equation}\label{eqninv}
 \pi (A) = \int_{\GS} \pi (dx) P(x, A),\quad A \in \cS.
\end{equation}

The transition probability defines the {\em transition operator} that acts by the formula 
\begin{equation}\label{troper}
 (Pf)(x) = \int_{\GS} P(x,dy) f(y)
\end{equation}
and is a positive contraction on every space $L^p(\pi) = L^p(\GS,\cS,\pi)$, $p\in [1,+\infty]$.

\subsection{$2$-U.I. condition}
 
Following \cite{Wu00} we will say that the transition operator $P$ is:
\begin{description}
\item{\bf uniformly integrable in $L^2$ (or $2$-U.I.)} if 
\begin{equation}\label{eq2ui}
\{ |Pf|^2\,;\, f\in L^2(\pi), \|f\|_2 \leq 1\} \ \ \text{ is uniformly $\pi$-integrable. }
\end{equation}
\item{\bf hyperbounded} if there exists $q > 2$ such that $P : L^2(\pi) \to
L^q(\pi)$ is a bounded linear operator, i.e.
\begin{equation}\label{eqhyp}
  \sup \{ \pi(|Pf|^q) \,;\, f\in L^2(\pi), \|f\|_2 \leq 1\} < +\infty.
\end{equation}
\item{\bf ultrabounded} if 
\begin{equation}\label{equltra}
 \sup \{ \|Pf\|_{\infty} \,;\, f\in L^1(\pi), \|f\|_1 \leq 1\} < +\infty.
\end{equation}
\end{description}

The hyperboundedness of the transition operator is, in a sense, independent of the particular choice of $p < q$, provided $1 < p < q < +\infty$. Indeed, by the Riesz-Thorin theorem, if $P$ is a bounded linear operator from $L^p$ to $L^q$, with $1 < p < q < +\infty$, then for any other $1 < p' < +\infty$ there is $q' > p'$, $q' < +\infty$,  such that $P$ is a bounded linear operator from $L^{p'}$ to $L^{q'}$. Notice also that if $P$ is ultrabounded, then for any $p > 1$
\[
\sup \{ \|Pf\|_{\infty} \,;\, f\in L^p(\pi), \|f\|_p \leq 1\} < +\infty.\]
In particular, the ultraboundedness implies the hyperboundedness and the latter implies the uniform integrability in $L^2$. 

Conditions like (\ref{eq2ui}) - (\ref{equltra}) are usually considered in the context of  hypercontractivity of Markov semigroups and all examples mentioned in \cite{Wu00} (as well as most of examples in \cite{CaM-A14}) are related to the continuous time Markov processes analysis. 

In the present paper we deal with discrete time Markov chains and show that also in this more elementary setting there are natural examples of Markov chains with contracting properties of the transition operator describable by relations (\ref{eq2ui}) - (\ref{equltra}).

For example, suppose that $P$ is given by a density $p(x,y)$ with respect 
to $\pi$, i.e.
\[ Pf(x) = \int_{\GS} \pi(dy) p(x,y) f(y).\]
Then $P$ is ultrabounded if $p(x,y)$ is a bounded function in $(x,y)$ (as in the main  model in \cite{JKO09}), and it is hyperbounded if $p(x,y) \in L^q(\pi \times \pi)$ for some $q > 2$ (see \cite[p. 480]{CaM-A14}). In Section \ref{SecEx1} we shall provide an example of a countable-space Markov chain with $P$ that is $2$-U.I. but not hyperbounded. 

\begin{remark} By the linearity of $P$, if any of conditions (\ref{eq2ui})- (\ref{equltra}) holds for {\em real-valued} functions $f$, then it is satisfied also for {\em complex-valued} functions $f$. 
\end{remark}

\subsection{$L^2$-spectral gap, geometric ergodicity and strong mixing}

The transition operator $P$ is said to have an {\em $L^2$-spectral gap} if there is a number $a < 1$ such that
\[
\sup\{ \|Pf\|_{L^2(\pi)}\, ; \,\int_{\GS}f(x)d\pi(x)=0,\,\|f\|_{L^2(\pi)}\leq 1\}\leq a.
\]
By iteration we obtain for $f \in L^2_0(\pi) = \{ f\in L^2(\pi)\,;\, \pi(f) = \int_{\GS} f(x) \pi(dx) = 0\}$
\begin{equation} \label{eq:gap}
\|P^n f\|_{L^2(\pi)} \leq a^n \|f\|_{L^2(\pi)},\quad n=1,2,\ldots .
\end{equation}
 This means that $\{X_n\}$ satisfies ``an $L^2$ norm condition'' of \cite{Ros71} and by Theorem 2, p. 217, {\em ibid.}, a central limit theorem with the standard normalization $\sqrt{n}$ holds  for the stationary sequence $\Psi(X_0), \Psi(X_1), \ldots$ whenever 
$\int \Psi(x) \pi(dx) = 0$ and $\int \Psi^2(x) \pi(dx) < +\infty$. (A proof of this limit theorem that is preferred nowadays can be found e.g.  in \cite{GoLi78}). 

For reversible, $\psi$-irreducible and aperiodic  Markov chains the spectral gap property is known to be equivalent to {\em geometric ergodicity}, i.e. existence of
 $ 0< \rho <1$ and $C : \GS \to \GR^+$ such that
\[\| P^n(x,\cdot) - \pi\|_{TV} \leq C(x) \rho^n,\quad \text{for $\pi$-a.e. $x\in\GS$},\]
where $\| \cdot \|_{TV}$ is the total variance distance (see \cite[Theorem 2.1]{RoRo97}). If $\{X_n\}$ is not reversible, then the spectral gap property implies  the geometric ergodicity (see \cite[Theorem 1.3]{KoMe12}), but there are Markov chains that are geometrically ergodic and do not have an $L^2$ spectral gap  (see \cite[Theorem 1.4]{KoMe12}). It is remarkable that the central limit theorem need not hold for such Markov chains (see \cite{Bra83}, \cite{Hagg05}, \cite{Hagg06}). 

Notice that if one is interested in a central limit theorem to hold for particular instantaneous function of the underlying Markov chain, then sufficient conditions weaker than the $L^2$ spectral gap are known (see e.g. \cite{MaWo00}). 

It is well known that the geometric ergodicity of a Markov chain is  equivalent (under natural conditions) to the exponential absolute regularity (see e.g. \cite[Theorem 21.19, p. 325]{Bra07II}), hence implies also {\em the strong mixing at geometric rate}. In this paper we shall use only the following consequence of the last property.

Let $\{X_j\}$ be strongly mixing at geometric rate. Then there exists a number $0 \leq \eta <1$ such that for any bounded measurable complex-valued function $\chi$ on $(\GS,\cS)$ 
\begin{equation}\label{eq:expo}
 \Big|\bE \Big(\chi(X_i) - \bE \big(\chi(X_i) \big) \Big)\Big( \overline{\chi(X_j) - \bE \big(\chi(X_j)\big)} \Big)\Big|  \leq 2\pi \eta^{|i-j|} \|\chi\|_{\infty}^2, \ i,j \in \GN.
 \end{equation}
See  \cite[Theorem 4.5, p. 125]{Bra07I}.

\subsection{Stable limits}

In the present paper the limiting distribution $\mu$ will be {\em stable}  with exponent $\alpha \in (0,2)$.  It is well-known (see e.g. \cite{SaTa94} or \cite{JaSh03}) that its characteristic function admits the L\'evy-Khintchine representation 
\begin{equation}\label{eq:f1}
\hat{\mu}(\theta)  = \exp\Big( i\theta a^h +   
\int \big( e^{i\theta x} - 1 - i \theta  x \Idn_{\{|x| \leq h\}}\big)\, \nu_{\alpha,c_+,c_-}(dx)\Big),
\end{equation}
where $c_+, c_-\geq 0$, $c_+ + c_- > 0$ and $a^h \in \GR^1$,  the L\'evy measure $\nu_{\alpha,c_+,c_-}$ has the density 
\[ p_{\alpha, c_+,c_-}(x) = \alpha \Big( c_+ x^{-(\alpha +1)}  \Idn_{\{x > 0\}} +  c_- |x|^{-(\alpha +1)}  \Idn_{\{x < 0\}}\Big),\]
and $h > 0$ is a fixed level of truncation. 
We will denote the stable distribution with characteristic function (\ref{eq:f1}) by $\delta_{a^h}*\text{$c_h$-Poiss}(\alpha, c_+, c_-)$.

In the main results of the paper we shall consider somewhat less general limits $\mu_{\alpha}$ with characteristic function of the form
\begin{equation}\label{strstab}
\hat{\mu}_{\alpha}(\theta)  = \begin{cases}
\exp\big( \int \big(e^{i\theta x} - 1\big) \nu_{\alpha, c_+, c_-}(dx)\big), &\alpha \in (0,1);\\
\exp\big(\int \big(e^{i\theta x} - 1\big) \nu_{1, c, c}(dx)\big), 
&\alpha =1;\\
\exp\big( \int \big(e^{i\theta x} - 1 -i\theta x\big) \nu_{\alpha, c_+, c_-}(dx)\big), &\alpha \in (1,2).
\end{cases}
\end{equation}
A reader familiar with the terminology would observe that completing the 
above list with probability laws of the form $\delta_a * \mu_1$, $a \not= 0$, we obtain all {\em strictly stable} laws on $\GR^1$ 

Notice that the integrals under the exponents in (\ref{eq:f1}) or (\ref{strstab}) can be evaluated, but obtained this way formulas are usually meaningless within the limit theory.

\section{Results}\label{SecResults}

Let $\{X_n\}$ be a Markov chain on the space $(\GS,\cS)$ with a stationary distribution $\pi$. Define 
$ \cF_n = \sigma\{ X_j\,;\, j \leq n\}.$
 
We will study distributional limits for suitably normalized and centered partial sums of the form
\[ S_n = \sum_{j=1}^{n} \Psi(X_j),\]
where $\Psi : (\GS,\cS) \to (\GR^1,\cB^1)$ is a measurable function. 

We will assume that the probability law $\pi\circ \Psi^{-1}$ belongs to the domain of attraction of $\mu_{\alpha}$, $0<\alpha<2$.
This means (see e.g. \cite[Theorem 1a, p. 313]{Fell70}) that 
\begin{equation}\label{domattr1}
 \pi \big( x\,;\, |\Psi(x)| > t\big) = t^{-\alpha} \ell(t),
\end{equation}
where $\ell(t)$ is a slowly varying function as $t\to\infty$, and
there exist the limits
\begin{equation}\label{domattr2}
 \lim_{t\to \infty} \frac{ \pi \big( x\,;\, \Psi(x) > t\big)}{ \pi \big( x\,;\, |\Psi(x)| > t\big)} = \frac{c_+}{c_+ + c_-},\ \ \lim_{t\to \infty} \frac{ \pi \big( x\,;\, \Psi(x) < - t\big)}{ \pi \big( x\,;\, |\Psi(x)| > t\big)} = \frac{c_-}{c_+ + c_-}.
\end{equation}

\begin{theorem}\label{TheMain}
Let $\{X_n\}$ be a Markov chain on the space $(\GS,\cS)$, with the 
transition operator $P$ and a stationary distribution $\pi$. We assume 
that $P$ has a spectral gap and satisfies the $2$-U.I. condition.

Let $\Psi : (\GS,\cS) \to (\GR^1,\cB^1)$ be such that  $\pi \circ \Psi^{-1}$ belongs to the domain of attraction of the stable distribution $\mu_{\alpha}$, $\alpha \in (0,2)$
(i.e. both (\ref{domattr1}) and (\ref{domattr2}) are fulfilled).
Let $B_n \to \infty$ satisfies 
\begin{equation}\label{eq:been}
 \frac{n}{B^{\alpha}_n} \ell(B_n) \to c_+ + c_-.
\end{equation}

\begin{description}
\item{\bf(i)} If $\alpha \in (0,1)$ or  $\alpha = 1$ and $c_+ = c_- = c$ then  
\[ \frac{\Psi(X_1) + \Psi(X_2) + \ldots + \Psi(X_{n})}{B_n} \indist \mu_{\alpha}.\]
\item{\bf(ii)} If $\alpha \in (1,2)$, then
\[ \frac{\sum_{j=1}^n \Psi(X_j) - \bE\big(\Psi(X_j)\vert \cF_{j-1}\big) }{B_n} \indist \mu_{\alpha}.\]

\end{description}
\end{theorem}
\begin{remarks}
\begin{enumerate}
\item Let us notice that in (ii) the tails of conditional expectations may {\em a priori}  influence the form of the limit. But they do not.
\item It is worth stressing that for $\alpha = 1$ we need only that {\em the limit} is symmetric and not $\pi\circ \Psi^{-1}$ itself.
\end{enumerate}
\end{remarks}

\begin{corollary}
In assumptions of Theorem \ref{TheMain}, if $\alpha \in (1,2)$ and 
\[ \bE \big( \Psi(X_1)\big| \cF_0\big) = 0. \]
i.e. $\Psi(X_1), \Psi(X_2), \ldots $ form a martingale difference sequence, then  
\[ \frac{\Psi(X_1) + \Psi(X_2) + \ldots + \Psi(X_{n})}{B_n} \indist \mu_{\alpha}.\]
\end{corollary}

There is another important case where we may get rid of centering by conditional expectations. As shown in  
\cite{JKO09} and \cite{CaM-A14} such situation takes place when we assume the hyperboundedness in place of condition $2$-U.I. We extend considerably the result of \cite{JKO09} by weakening the $L^2$-spectral gap property to the strong mixing at geometric rate.

\begin{theorem}\label{TheHyper}
In assumptions of Theorem \ref{TheMain} replace the $2$-U.I. condition with the hyperboundedness and the $L^2$-spectral gap property with the strong mixing at geometric rate (in particular: with the geometric ergodicity).

Then 
\[ \frac{\Psi(X_1) + \Psi(X_2) + \ldots + \Psi(X_{n})}{B_n} \indist \mu_{\alpha},\]
provided 
\begin{description}
	\item{\bf(i)} $\alpha \in (0,1)$; 
	\item{\bf(ii)} $\alpha = 1$ and $c_+ = c_- = c$, 
	\item{\bf(iii)} $\alpha \in (1,2)$ and $\int \Psi(x) \pi (dx) = 0$. 
\end{description}
\end{theorem}

As a by-product of developed techniques we obtain a weak law of large numbers for geometrically ergodic Markov chains, which might be of independent interest.

\begin{theorem}\label{TheWeak}
Let $\{X_n\}$ be strongly mixing at geometric rate (in particular: geometrically ergodic) Markov chain on $(\GS,\cS)$ with a stationary distribution $\pi$. Suppose  $\Psi : (\GS,\cS) \to (\GR^1, \cB^1)$ is such that for some $\beta > 1$ 
\[ \int \pi(dx) |\Psi(x)|^{\beta} < +\infty,\quad \int \pi(dx) \Psi(x) = 0.\] 
Then for any $\alpha \in (0,\beta \wedge 2)$ and any $1/\alpha$-regularly varying sequence $B_n$ we have
\[ \frac{\Psi(X_0) + \Psi(X_1) + \ldots + \Psi(X_{n-1})}{B_n} \inprob 0.\] 
\end{theorem}

\section{Examples}\label{SecExamples}

\subsection{Example related to the $2$-U.I. condition}\label{SecEx1}

We are going to construct a discrete in time and space example of the transition operator that exhibits the $L^2$-spectral gap property, satisfies the $2$-U.I. condition but is not hyperbounded. This will show that our theory essentially extends the results of \cite{JKO09} and \cite{CaM-A14}. Notice also that all the examples of operators provided in \cite{Wu00} and related to the $2$-U.I. condition are taken from the stochastic analysis.

\begin{example}\label{examain}     
The example is a variant of  Rosenblatt's family of examples  \cite[ pp. 213-214]{Ros71}, but it occurs also in many other places, e.g. in \cite[p. 54]{MeTw09}, in the context of the backward recurrence time chain.

Let $T: (\Omega,\cF,\bP) \to \GN = \{0,1,2,\ldots\}$ be an integer valued nonegative random variable such that 
\[ \bE T < +\infty,\ \bP(T \geq j) > 0,\ j\in \GN.\]
(Other requirements imposed on the distribution of $T$ will be specified later). Let the transition probabilities $p_{j,k}$ be given by the formula
\[ p_{j,k} = \begin{cases}
\displaystyle{\frac{\bP(T=j)}{\bP(T \geq j)},} &\text{ if $k=0$};\\
\displaystyle{\frac{\bP(T\geq j+1)}{\bP(T \geq j)}}, &\text{ if $k=j+1$};\\
0, &\text{ otherwise}.
\end{cases}\]
Then 
\[ \pi(j) = \frac{\bP(T \geq j)}{1 + \bE T},\ j=0,1,2, \ldots,\]
is the unique stationary distribution for $P = [p_{j,k}]$ and the transition operator reads
\[ Pf(j) = \frac{\bP(T=j)}{\bP(T \geq j)} f(0) + \frac{\bP(T\geq j+1)}{\bP(T \geq j)} f(j+1). \] 
Let $\{X_n\}$ be a Markov chain on $\GS = \GN$  with the transition probabilities $[p_{j,k}]$.
\begin{lemma}\label{lem:gap}
If\[ 3\, \bE T < \bP ( T = 0),\]
and
\begin{equation}\label{eq:domin}
 \bP( T\geq 1) \geq \sup_{k \geq 1} \frac{\bP( T \geq k+1)}{\bP( T \geq k)},
 \end{equation}
then the Markov chain $\{X_n\}$ has the $L^2$-spectral gap property. 
\end{lemma}
\noindent {\sc Proof.}\ Let $f \in L^2_0(\pi)$ and $\| f\|_2 =1$. These relations imply that
\begin{align*} \big|f(0)\big| &= \Big|- \sum_{j=1}^{\infty} f(j)\bP( T\geq j)\Big| \leq \sum_{j=1}^{\infty} \big(|f(j)|\sqrt{\bP( T\geq j)}\big)
\sqrt{\bP( T\geq j)} \\
&\leq \sqrt{\sum_{j=1}^{\infty} f^2(j) \bP( T \geq j)} \sqrt{\sum_{j=1}^{\infty} \bP( T \geq j)}  = \sqrt{\big(1 + \bE T - f^2(0)\big) \bE T}. 
\end{align*}
Hence 
\begin{equation}\label{eq:efzero}
|f(0)| \leq \sqrt{\bE T}.
\end{equation}
In a similar way we obtain
\begin{equation}\label{eq:sumefzero}\begin{split}
\Big|\sum_{j=1}^{\infty} f(j) &\frac{\bP (T=j-1)}{\bP(T \geq j-1)} \frac{\bP(T \geq j)}{1 + \bE T}\Big| \leq \sum_{j=1}^{\infty} |f(j)|  \frac{\bP(T \geq j)}{1 + \bE T} \\
 & \leq \sqrt{\sum_{j=1}^{\infty} f^2(j) \frac{\bP( T \geq j)}{1 + \bE T}} \sqrt{\sum_{j=1}^{\infty} \frac{\bP( T \geq j)}{1 + \bE T}} \leq \sqrt{\frac{\bE T}{1 + \bE T}}.
\end{split}
\end{equation}
We are ready for estimates of $\mathbb{E}_{\pi}\big(\vert Pf\vert^2\big)=\big(1/(1 + \bE T)) \sum_{j=0}^{\infty} \big| Pf (j)\big|^2 \bP( T\geq j)$.
\begin{align*}
\sum_{j=0}^{\infty} & \big| Pf (j)\big|^2 \frac{\bP( T\geq j)}{1 + \bE T} = \sum_{j=0}^{\infty} \Big|\frac{\bP(T=j)}{\bP(T \geq j)} f(0) + \frac{\bP(T\geq j+1)}{\bP(T \geq j)} f(j+1)\Big|^2 \frac{\bP( T\geq j)}{1 + \bE T} \\
& = \frac{f^2(0)}{1 + \bE T} \sum_{j=0}^{\infty} \frac{ \bP^2(T=j)}{\bP(T\geq j)} + \frac{2 f(0)}{1 + \bE T} \sum_{j=0}^{\infty} f(j+1) \frac{\bP( T = j)}{\bP(T\geq j)} \bP( T \geq j+1) \\
&\qquad + \frac{1}{1 + \bE T} \sum_{j=0}^{\infty} f^2(j+1) \frac{ \bP^2 ( T \geq j+1)}{\bP( T \geq j)} = J_1 + J_2 + J_3.
\end{align*}
We have by (\ref{eq:efzero}) 
\[ J_1 \leq \frac{\bE T}{1 + \bE T} \sum_{j=0}^{\infty} \frac{ \bP^2(T=j)}{\bP(T\geq j)}  \leq \bE T,\] 
while by  (\ref{eq:efzero}) and (\ref{eq:sumefzero})
\[ J_2 \leq 2 \sqrt{\bE T} \sqrt{\frac{\bE T}{1 + \bE T}} \leq 2 \bE T. \]
Finally, by (\ref{eq:domin}),
\[ J_3 \leq \bP (T \geq 1) \sum_{j=1}^{\infty} f^2(j) \frac{\bP(T \geq j)}{1 + \bE T} \leq \bP( T \geq 1).   \]
Therefore 
\[ \mathbb E_{\pi}\left(\vert Pf\vert^2\right)\leq 3\, \bE T + \bP( T \geq 1) = 1 - \big(\bP( T =0) - 3 \bE T\big)  = a < 1. \]
The proof of Lemma 4.2 is complete.

It remains to show that for some specific distribution of $T$ the $2$-U.I. condition holds, but there is no hyperboundedness.
Choose $\gamma \in (0,1)$ and set
\[ \bP( T\geq 1) = \gamma, \bP( T\geq 2) = \gamma^3, \ldots, \bP(T \geq  j) = \gamma^{1 + 2 + \ldots + j} = \gamma^{j(j+1)/2},\ldots.\]
Clearly, $\bP( T \geq j+1)/\bP( T \geq j) = \gamma^{j+1},\ j=0,1,2,\ldots$ and for $\gamma < 1/5$
\[\bE T < \frac{\gamma}{1 - \gamma} < (1/3)(1 - \gamma) = (1/3) \bP( T =0),\]
so that the assumptions of Lemma \ref{lem:gap} are satisfied and the corresponding Markov chain $\{X_n\}$ has the $L^2$-spectral gap property. 

In order to prove that the $2$-U.I. condition holds, it is enough to show that
\[ \lim_{k\to\infty} \sup_{\|f\|_2 \leq 1}\sum_{j=k}^{\infty}  \big| Pf (j)\big|^2 \frac{\bP( T\geq j)}{1 + \bE T} = 0.\]
Notice that $\|f\|_2 \leq 1$ implies that $f^2(j) \leq (1 + \bE T)/ \bP(T \geq j), j = 0, 1, 2, \ldots $ Keeping this in mind we can proceed as follows.
\begin{align*}
\sum_{j=k}^{\infty} \big| Pf (j)\big|^2 &\frac{\bP( T\geq j)}{1 + \bE T} \\
& \leq \frac{2 f^2(0)}{1 + \bE T} \sum_{j=k}^{\infty} \frac{ \bP^2(T=j)}{\bP(T\geq j)} + \frac{2}{1 + \bE T}  \sum_{j=k}^{\infty} f^2(j+1) \frac{ \bP^2 ( T \geq j+1)}{\bP( T \geq j)} \\
&\leq 2 \bP(T \geq k) + 2 \sum_{j=k}^{\infty}  \frac{ \bP ( T \geq j+1)}{\bP( T \geq j)} =  2 \bP(T \geq k) + 2 \sum_{j=k}^{\infty} \gamma^{j+1} \to 0.
\end{align*}

Next consider a sequence $\{f_k\}$ of functions in $L^2(\pi)$ given by
\[ f_k(j) = \begin{cases}
\displaystyle{\sqrt{\frac{1 + \bE T}{\bP( T \geq k)}}}, &\text{ if $j=k$};\\
0, &\text{ otherwise}.
\end{cases}\]
Take any $q > 2$. We have, if $k \to \infty$,
\begin{align*}
 \|P f_k\|^q_q &= \sum_{j=0}^{\infty} \big| Pf_k(j)\big|^q \frac{\bP( T\geq j)}{1 + \bE T} \\
& =  \Bigg( \frac{1 + \bE T}{\bP( T \geq k)}\Bigg)^{q/2} \Bigg(\frac{\bP( T \geq k)}{\bP(T \geq k-1)}\Bigg)^q \frac{\bP( T \geq k-1)}{1 + \bE T} \\
& = \big( 1 + \bE T\big)^{q/2 - 1}\frac{ \big(\bP(T \geq k)\big)^{q/2}}{ \big(\bP(T \geq k-1)\big)^{q -1}} = \big( 1 + \bE T\big)^{q/2 - 1} \gamma^{w(k)} \to +\infty,
\end{align*}
for $w(k) = q k(k+1)/4 - (q - 1) k(k-1)/2 = (1/4) \big( k^2(2-q) + k (3q - 2)\big) \to -\infty$. It follows that the transition operator $P$ cannot be a bounded linear map from $L^2(\pi)$ to $L^q(\pi)$. 
\end{example}

\subsection{Gaussian hyperboundedness}
Let us examine a standard example, already considered by Doob \cite[p.218]{Doob53}. For $0 < |\rho| < 1$ set 
\[ P(x, dy) = \frac{1}{\sqrt{2\pi(1 - \rho^2)}} e^{-\frac{(y - \rho x)^2}{2(1-\rho^2)}} \,dy,\]
and notice that for each $x\in\GR^1$
\[ P^n(x, dy) = \frac{1}{\sqrt{2\pi(1 - \rho^{2n})}} e^{-\frac{(y - \rho^n x)^2}{2(1-\rho^{2n})}} \,dy \conver_{TV} \pi(dy) = \frac{1}{\sqrt{2\pi}} e^{-\frac{y^2}{2}} \,dy. \]
Rosenblatt \cite[p. 214]{Ros71} proves that the Markov chain $\{X_n\}$ corresponding to $P(x,dy)$ has the $L^2$-spectral gap property. Alternatively, we may observe that $\{X_n\}$ is a Gaussian stationary sequence with the correlation function $\bE X_0 X_k = \rho^{|k|}$, hence an $AR(1)$ process with Gaussian innovations. By \cite[p.389]{MeTw09} it is a geometrically ergodic Markov chain. Since it is also reversible, it admits an $L^2$-spectral gap. 

Using \cite[p. 480]{CaM-A14}), we shall show that $\{X_n\}$ is also hyperbounded. Indeed, $P(x, dy) = p(x,y)\pi(dy)$, where
\[ p(x,y) = \frac{1}{\sqrt{1 - \rho^2}} \exp\Big( -\frac{\rho^2}{2(1 - \rho^2)}x^2 + \frac{\rho x y}{1-\rho^2} -
\frac{\rho^2}{2(1 - \rho^2)}y^2\Big).\] 
And we have 
\[ \int \pi(dx) \pi(dy) p(x,y)^q < +\infty,\]
whenever 
\[2 < q < \frac{1 + |\rho|}{|\rho|}.\]
Hence we can apply Theorem \ref{TheHyper} and obtain 
\[ \frac{\Psi(X_1) + \Psi(X_2) + \ldots + \Psi(X_{n})}{B_n} \indist \mu_{\alpha},\]
for suitably chosen $\Psi$ and $B_n$. Notice that the lack of any necessary centering confirms in this particular case the conjecture of Davis formulated in the frame of Example on p. 267 in \cite{Dav83}.

There is another reason for evoking this classic example.
It was Doob \cite[p.218]{Doob53} who pointed out that this Markov chain does not satisfy Doeblin's condition (D). And since the work of Davydov \cite{Davy73} we know that Doeblin's condition means essentially $\phi$-mixing of a Markov chain. 
It follows that the limit theory developed in our paper is much broader than results depending on uniform ergodicity of Markov chains, as presented e.g. in \cite{Cog72}. 

\subsection{ARCH processes with heavy tails are not hyperbounded}

An ARCH($1$) process is a Markov chain given by the recurrence
formula
\begin{equation}\label{ajarch}
X_{j+1} = \sqrt{\beta + \lambda X_j^2} Z_{j+1}, \   j \geq 0,
\end{equation}
where $\beta,\lambda > 0$ and $\{Z_n\}_{n\in\GN}$ is 
an i.i.d. sequence, independent of $X_0$. In order 
to comply with references we shall assume that 
$Z_n \sim \cN(0,1)$. 

For basic information on ARCH processes and the properties used below  we refer both to the classic book \cite{EKM97} and to the recent source \cite{BDM16}.

In the range of parameters $\beta > 0$
and $\lambda \in (0,2 e^{\gamma})$ (where $\gamma$ is the Euler constant) the process 
$\{X_j\}$ admits a stationary distribution given by 
\[ X_0 \sim r_0 \sqrt{\beta \sum_{m=1}^{\infty} 
	Z_m^2 \prod_{j=1}^{m-1}
	(\lambda Z_j^2)},\]
where $r_0$ is a Rademacher random variable 
($P(r_0 = \pm 1) = 1/2$), independent of $\{Z_n\}$.
This stationary distribution exhibits power decay of the tails. Namely, if $\kappa > 0$ is
the unique positive solution of the equation
\[ \bE (\lambda Z_1^2)^u = 1,\]
then, as $x \to \infty$,
\begin{equation}\label{aj11}
\bP\big( X_0 > x\big) = \bP\big( X_0 < -x\big)\sim \frac{C_{\beta,\lambda}}2\ x^{-2\kappa},
\end{equation}
where
\[ C_{\beta,\lambda} = \frac{\bE\Big[\big(\beta + \lambda X^2_0\big)^{\kappa} -
	\big(\lambda X^2_0\big)^{\kappa}
	\Big]}{\kappa \lambda^{2\kappa} E\Big[\big(Z_1^{2\kappa}
	\ln (\lambda Z_1^2)\Big]} \in (0,+\infty).\]
It follows that $\lambda > 1$ implies ``really'' heavy tails and it is likely that the partial sums of $\{X_j\}$ properly normalized converge to stable laws. Indeed, Davis and Mikosch \cite{DaMi98} showed that the  partial sums under the natural normalization  converge to {\em some} stable limit and Bartkiewicz et al. \cite{BJMW11} identified the parameters of the limit. 

For purposes of the present example, let us denote by
$\mu_{\alpha,\tau}$ the 
symmetric
$\alpha$-stable distribution given for $\alpha \in (0,2)$ and 
$\tau > 0$ by
\[\widehat{\mu_{\alpha,\tau}}(\theta) =
\exp\Big(\tau \alpha \int_{\GR}\big(e^{i\theta u} -
1\big)|u|^{-(\alpha+1)}\,du\Big). \] 
If our Theorems \ref{TheMain} or \ref{TheHyper} were applicable to $\{X_j\}_{j\geq 0}$, then we would have 
\[\frac{X_1 + X_2 + \ldots + X_{n}}{( n 
	C_{\beta,\lambda})^{\frac{1}{2\kappa}}}
\indist \mu_{2\kappa,1}.\]
It is, however,  proved in \cite{BJMW11} that 
	\[\frac{X_1 + X_2 + \ldots + X_{n}}{( n C_{\beta,\lambda})^{\frac{1}{2\kappa}}}
\indist \mu_{2\kappa,\tau},\]
where $\tau = E \big[|1 +
S_{\infty}|^{2\kappa} - |S_{\infty}|^{2\kappa}\big] > 0$
and the series
\[S_{\infty} =
\sum_{j=1}^{\infty} \lambda^{j/2}\big[\prod_{k=1}^{j-1}|Z_k|\big]Z_j\]
converges a.s.

Moreover, ARCH($1$) processes are strongly mixing at geometric rate, as is shown in \cite[p. 2077]{DaMi98}.

Therefore the transition operator of an ARCH($1$) process {\em is not hyperbounded}.

\begin{remark}
	If $2\kappa \in (1,2)$, then the corresponding ARCH($1$) process $\{X_j\}$ forms a martingale difference sequence, partial sums of which normalized by $n^{1/{2\kappa}}$ are weakly convergent, but to a different limit than in the independent case. This is in striking contrast to the properties of martingale difference sequences with finite variance!
\end{remark}

\subsection{$m$-skeletons}
It is well known that iterating the transition operator improves its properties from many viewpoints. 
So it may happen that some power  $P^m$ is hyperbounded, 
for instance, while $P$ itself not. Such situation implies 
that for $\{\Psi(X_{k\cdot m})\}_{k=0,1,2,,\ldots}$ (the $m$-skeleton)
 some $\alpha$-stable limit theorem holds and one may hope to 
 extend this property to the whole sequence. This is 
 impossible in general,  as the simple counterexample 
 provided already by Rosenblatt
  \cite[p. 195]{Ros71} shows. Indeed, take an i.i.d. 
  sequence $\{Y_n\}$ of strictly stable random variables 
  and consider a Markov chain on $\GS = \GR^2$ given by 
  the formula $X_n = (Y_n, Y_{n-1})$. Take 
  $\Psi(x,y) = x-y$. Then $\sum_{j=0}^{n-1} \Psi(X_n)$ 
  remains stochastically bounded while the $1$-skeleton
   consists of independent random variables and therefore 
   satisfies the corresponding limit theorem.
 
 Rosenblatt's example is of probabilistic provenience. Some people may prefer another example given below that is closer to thinking in terms of dynamical systems.
 
 \begin{example}
Set $\GS = [0,3)$ and let $Leb$ be the Lebesgue measure restricted to $\GS$. For $x\in [0,1)$ and  $B \in \cB_{[0,1)} \cup \cB_{[2,3)}$ define
\[ P(x, \{x+1\}) = P(x+1, \{x+2\})  = 1,\ P(x+2, B) = \frac{Leb(B)}{2}.\]
The invariant measure $\pi$ is given by the density
\[ p(x) = \frac{1}{4}\Idn_{[0,2)}(x)  +  \frac{1}{2}\Idn_{[2,3)}(x).\]
Elementary calculations show that for $f \in L^2_0(\pi)$ we have
\[ \bE_{\pi} \Big( \big(P^3f\big)^2\Big) \leq \frac{27}{32} \bE_{\pi}\big(f^2\big),\]
i.e. the $3$-skeleton has the spectral gap property. Another elementary calculation shows that also
\[ \|P^3f\|_{\infty} \leq  3\|f\|_1,\]
i.e. the $3$-skeleton is ultrabounded. 

Now take $\psi(\cdot) : [0,1) \to \GR^1$ with a symmetric $\alpha$-stable distribution $\mu$ and define 
\[ \Psi(x) = \begin{cases}
\psi(x), &\text{ if $x \in [0,1)$};\\
- \psi(x-1), &\text{ if $x \in [1,2)$};\\
0, &\text{ if $x \in [2,3)$}.
\end{cases}
\]
One  verifies directly that
\[ \pi\big( \Psi > r\big) = \frac{1}{2} Leb\big( \psi > r\big),
 \quad \pi\big( \Psi < - r\big) = \frac{1}{2} 
 Leb\big( \psi < - r\big).\]
Therefore the $3$-skeleton $\{\Psi(X_{k3})\}$ satisfies all
 assumptions of our Theorem \ref{TheHyper}, while the 
 partial sums of the whole sequence $\{\Psi(X_k)\}$ are 
 bounded in probability. 
 \end{example}

\section{Proofs}\label{SecProofs}
\subsection{Some auxiliary results}
We begin with establishing an important property of conditional distributions $P(x,\,dy) \circ \Psi^{-1}$ that is a consequence of solely (\ref{domattr1})-(\ref{domattr2}).

\begin{proposition}\label{PropFirst} Suppose that (\ref{domattr1}) and (\ref{domattr2}) hold. Let $B_n$ be defined by (\ref{eq:been}). Then
\begin{equation}\label{eq:sqinp}
n\Big| 1 - \bE\big( e^{i\theta \Psi(X_1)/B_n}\big| \cF_{0}\big)\Big|^2 \inprob 0,\quad \theta \in \GR^1.
\end{equation}
\end{proposition}
\noindent{\sc Proof.}\ Recall that if $B_n$ is defined by (\ref{eq:been})
 then $B_n = n^{1/\alpha} \tilde{\ell}(n)$, where $\tilde{\ell}(t)$ is a 
slowly varying function. Let $h > 0$ be fixed. Using the inequality  $| 1 + ix - e^{ix}| \leq \frac12 |x|^2$, we have
\begin{align}
n\Big| 1& - \bE\big( e^{i\theta \Psi(X_1)/B_n}\big| \cF_{0}\big)\Big|^2 \leq \nonumber \\
&\leq 2 n\Big| 1 +  i\theta\bE \big(\frac{\Psi(X_1)}{B_n} \Idn_{\{|\Psi(X_1)|
	 \leq h B_n\}}\big| \cF_{0}\big) - \bE\big( e^{i\theta \Psi(X_1)/B_n}\big| \cF_{0}\big)\Big|^2 \nonumber\\
&\quad\quad + 2 n B_n^{-2} \theta^2\Big|\bE \big(\Psi(X_1) \Idn_{\{|\Psi(X_1)| \leq h B_n\}} \big|
 \cF_{0}\big)\Big|^2 \label{eq:basic}\\
&\leq n B_n^{-4} \theta^4 \Big(\bE \big(\Psi(X_1)^2 \Idn_{\{|\Psi(X_1)| \leq h B_n\}} \big| \cF_{0}\Big)^2 + 16 n \Big( \bP \big( |\Psi(X_1)| > h B_n \big| \cF_0\big)\Big)^2 + \nonumber\\
&\quad\quad +  2 n B_n^{-2} \theta^2\Big|\bE \big(\Psi(X_1) \Idn_{\{|\Psi(X_1)| \leq h B_n\}} \big| \cF_{0}\big)\Big|^2
= \theta^4 I^h_{n,1} + 16 I^h_{n,2} + 2\theta^2 I^h_{n,3}.\nonumber
\end{align}
At first we shall examine the convergence of $I^h_{n,3}$. 
If $\alpha \in (1,2)$ then 
\[ \Big|\bE \big(\Psi(X_1) \Idn_{\{|\Psi(X_1)| \leq h B_n\}} \big| \cF_{0}\big)\Big|^2 \to \Big|\bE \big(\Psi(X_1) \big| \cF_{0}\big)\Big|^2\ \text{ a.s., }\]
while $n B_n^{-2} = n^{1 - 2/\alpha}(\tilde{\ell}(n))^{-2} \to 0$. Consequently,  $I^h_{n,3} \to 0$ a.s.

Now suppose that $\alpha \in (0,1]$. 
Take  $0 < r <\alpha/2$. We have 
\[ \bE \Big(\bE\big( |\Psi(X_1)|^{\alpha - r}\big| \cF_0\big) \Big) = \bE |\Psi(X_1)|^{\alpha - r} < +\infty,\]
and so 
\[ (\alpha -r) \int_{0}^{\infty} t^{\alpha -r -1} P\big(X_0, |\Psi|^{-1}(t,+\infty) \big) \,dt = \bE\big( |\Psi(X_1)|^{\alpha - r}\big| \cF_0\big) < +\infty\ \text{a.s.}\]
It follows that
\begin{align*}
I^h_{n,3} &= n B_n^{-2} \Big|\bE \big(\Psi(X_1) \Idn_{\{|\Psi(X_1)| \leq h B_n\}} \big| \cF_{0}\big)\Big|^2 \\
&\leq n B_n^{-2} \Big|\bE \big(\big|\Psi(X_1)\big| \Idn_{\{|\Psi(X_1)| \leq h B_n\}} \big| \cF_{0}\big)\Big|^2\\
&\leq  n B_n^{-2} \Big| \int_{0}^{h B_n} P\big(X_0, |\Psi|^{-1}(t,+\infty) \big) \,dt \Big|^2\\
 &=  n B_n^{-2} \Big| \int_{0}^{h B_n} t^{1 - \alpha + r} t^{\alpha - r - 1} P\big(X_0, |\Psi|^{-1}(t,+\infty) \big) \,dt \Big|^2 \\
&\leq n B_n^{-2}  h^{2(1 - \alpha +r)} B_n^{2(1 - \alpha +r)} \Big| \int_{0}^{\infty} t^{\alpha - r - 1}P\big(X_0, |\Psi|^{-1}(t,+\infty) \big) \,dt \Big|^2 
\\
& = n^{-1 + 2r/\alpha} \big(\tilde{\ell}(n)\big)^{-2(\alpha-r)}  h^{2(1 - \alpha +r)}\Big(\frac{1}{\alpha - r} \bE\big( |\Psi(X_1)|^{\alpha - r}\big| \cF_0\big) \Big)^2  \to 0, \ \text{a.s.}
\end{align*}
Similarly, if $\alpha \in (0,2)$ and $0 < r <\alpha/2$, then we have 
\begin{align*}
I^h_{n,1} &= n B_n^{-4} \Big|\bE \big(\Psi(X_1)^2 \Idn_{\{|\Psi(X_1)| \leq h B_n\}} \big| \cF_{0}\big)\Big|^2 \\
& \leq 4 n B_n^{-4} \Big(\int_{0}^{h B_n} t P\big(X_0, |\Psi|^{-1}(t,+\infty) \big) \,dt \Big|^2\\
&= 4 n B_n^{-4} \Big(\int_{0}^{h B_n} t^{2 - \alpha + r} t^{\alpha - r - 1} P\big(X_0, |\Psi|^{-1}(t,+\infty) \big) \,dt \Big)^2 \\
&\leq 4 n B_n^{-4} h^{2(2 - \alpha + r)} B_n^{2(2 - \alpha + r)}
\Big(\int_{0}^{\infty} t^{\alpha - r - 1} P\big(X_0, |\Psi|^{-1}(t,+\infty) \big) \,dt \Big)^2 \\
&= 4 n^{-1 + 2r/\alpha} \big(\tilde{\ell}(n)\big)^{-2(\alpha-r)}  h^{2(2 - \alpha +r)}\Big(\frac{1}{\alpha - r} \bE\big( |\Psi(X_1)|^{\alpha - r}\big| \cF_0\big) \Big)^2 \to 0, \ \text{a.s.}
\end{align*}
It remains to show that  $I^h_{n,2} \inprob 0$. This condition is not related to truncated moments and therefore requires a different type argument. Notice that the convergence in probability is metrizable and so it is enough to show that in every subsequence $n'$ one can find a further subsequence $n^{\prime\prime}$ along which $I^h_{n^{\prime\prime},2} \inprob 0$. So choose $n'$ and consider  random variables $Y_{n'}$ defined on $(\GS,\cS)$ by the formula 
\[ Y_{n'}(x) = n' P\big(x, |\Psi|^{-1}(h B_{n'},+\infty)\big).\]
We know from (\ref{domattr1}), (\ref{domattr2}), (\ref{eq:been}) and the continuity of the stable L\'evy measure that 
\[ \int_{\GS} \pi(dx)   Y_{n'}(x) = n'\bP\big(|\Psi(X_1)| > h B_{n'}\big) \to 
(c_+ + c_-) h^{-\alpha},\]
hence, in particular, random variables $\{Y_{n'}\}$ are uniformly tight. 
Let $\{n^{\prime\prime}\}$ be a subsequence such that 
$Y_{n^{\prime\prime}} \indist Y_{\infty}$. By the Skorokhod 
representation theorem one can construct random variables 
$\widetilde{Y}_{n^{\prime\prime}}$ and $\widetilde{Y}_{\infty}$,
 defined on the standard probability space
$\big([0,1], \cB_{[0,1]}, Leb\big)$ and such that
\[ \widetilde{Y}_{n^{\prime\prime}} \sim Y_{n^{\prime\prime}},\quad  \widetilde{Y}_{\infty} \sim Y_{\infty},\]
and
\[ \widetilde{Y}_{n^{\prime\prime}}(\omega) \to \widetilde{Y}_{\infty}(\omega),\quad \text{for almost all $\omega \in [0,1]$.}\]
This implies that
\[ \frac{1}{n^{\prime\prime}}\widetilde{Y}^2_{n^{\prime\prime}}(\omega) \to 0,\quad \text{for almost all $\omega \in [0,1]$.}\]
 But under the initial distribution $\pi$ we have
\begin{equation}\label{eq:rate}
n^{\prime\prime} \Big(P\big(x, |\Psi|^{-1}(h B_{n^{\prime\prime}},+\infty)\big)\Big)^2 = \frac{1}{n^{\prime\prime}}Y^2_{n^{\prime\prime}} \sim \frac{1}{n^{\prime\prime}}\widetilde{Y}^2_{n^{\prime\prime}}.
\end{equation}
It follows that 
\[  n^{\prime\prime} \Big(P\big(x, |\Psi|^{-1}(h B_{n^{\prime\prime}},+\infty)\big)\Big)^2 \inprob 0.\]
\begin{remark}
It is clear that the convergences $I^h_{n,1} \inprob 0$ and $I^h_{n,3} \inprob 0$ can be obtained also by the last method. But the proofs given above lead to the a.s. convergence and provide some idea about the rate of convergence.
\end{remark}
\begin{remark}\label{rem:rate}
	It is also clear that relation (\ref{eq:rate}) can be extended to 
\[	(n^{\prime\prime})^{1+\delta} \Big(P\big(x, |\Psi|^{-1}(h B_{n^{\prime\prime}},+\infty)\big)\Big)^2 = \frac{1}{(n^{\prime\prime})^{1 - \delta}}Y^2_{n^{\prime\prime}} \sim \frac{1}{(n^{\prime\prime})^{1 - \delta}}\widetilde{Y}^2_{n^{\prime\prime}},\]
hence, in fact, we have 
\[ n^{\delta} I^h_{n,2} \inprob 0, \]
for every $\delta \in [0,1)$. Gathering information on $I_{n,1}^h, I_{n,2}^h$ and $I_{n,3}^h$ we obtain existence of some $\delta > 0$ such that 
\[n^{1 + \delta} \Big| 1 - \bE\big( e^{i\theta \Psi(X_1)/B_n}\big| \cF_{0}\big)\Big|^2 \inprob 0,\quad \theta \in \GR^1.\]
	\end{remark}

Now we are ready to prove two universal (i.e. independent of $\alpha \in (0,2)$) limit theorems.

\begin{proposition}\label{PropMain}
Let $\{X_n\}$ be a Markov chain on the space $(\GS,\cS)$, with the 
transition operator $P$ and a stationary distribution $\pi$. We assume 
that $P$ has a spectral gap and satisfies the $2$-U.I. condition.

Let $\alpha \in (0,2)$ and $h > 0$. Let $\Psi : (\GS,\cS) \to (\GR^1,\cB^1)$ be such that  $\pi \circ \Psi^{-1}$ belongs to the domain of attraction of the stable distribution $\mu_{\alpha}$, $\alpha \in (0,2)$
(i.e. both (\ref{domattr1}) and (\ref{domattr2}) are fulfilled). 
Let $B_n \to \infty$ satisfies 
\[ \frac{n}{B^{\alpha}_n} \ell(B_n) \to c_+ + c_-.\]
Set $S^h_n =  \sum_{j=1}^{n} \Psi(X_j) - \bE \big(\Psi(X_j) \Idnem_{\{|\Psi(X_j)| \leq h B_n\}} \big| \cF_{j-1}\big)$.
Then
\begin{equation}\label{eq:cent}
\frac{S^h_n}{B_n} \indist \text{$c_h$-Poiss}(\alpha, c_+, c_-).
\end{equation}
\end{proposition}

\noindent{\sc Proof.}\  Choose $\theta \in \GR^1$ and 
notice that by Proposition \ref{PropFirst} relation  (\ref{eq:sqinp}) holds. 
We will show that this relation can be strengthened to 
\begin{equation}\label{eq:sqinelone}
n \bE \Big| 1 - \bE\big( e^{i\theta \Psi(X_1)/B_n}\big| \cF_{0}\big)\Big|^2 
\to 0.
\end{equation}
It is enough to  show that $ n \big| 1 - \bE\big( e^{i\theta \Psi(X_1)/B_n}\big| \cF_{0}\big)\big|^2$ is a uniformly integrable sequence. By the $2$-U.I. condition we have to prove that the sequence $\{Z_n = \sqrt{n}  \big(1 - e^{i\theta \Psi(X_1)/B_n}\big)\}$ is bounded in $L^2$.
\begin{align}
\bE \Big| \sqrt{n}  &\big(1 - e^{i\theta \Psi(X_1)/B_n}\big) \Big|^2 \nonumber\\
 &=
n \,\bE  \Big(\big(1 - \cos(\theta \Psi(X_1)/B_n)\big)^2 +  \big(\sin(\theta \Psi(X_1)/B_n)\big)^2\Big) \label{eq:bound} \\
&\leq \theta^2(1 + \theta^2/4) \frac{n}{B_n^2} \bE \Psi(X_1)^2 \Idn_{\{|\Psi(X_1)| \leq B_n\}}  + 5 n \bP\big( |\Psi(X_1)| > B_n \big) \nonumber \\
&\leq \theta^2(1 + \theta^2/4) \frac{n}{B_n^2} 2 \int_0^{B_n} t \bP\big( |\Psi(X_1)| > t\big)\,dt + 5 n \bP\big( |\Psi(X_1)| > B_n \big). \nonumber
\end{align}
 The last  expression converges to $\big(2\theta^2(1 + \theta^2/4) + 5\big)(c_+ + c_-)$ by the definition of $B_n$ and the direct half of the Karamata theorem (see \cite[Theorem 1.5.11, p. 28]{BGT87}). 

Given (\ref{eq:sqinelone}) we obtain the crucial relation (\ref{eq:new1}) 
\[ \bE\Big(\sum_{j=1}^{n} \Big| 1 - \bE\left(e^{i\theta \Psi(X_j)/B_n}\vert\cF_{j-1}\right) \Big|^2\Big) = n \bE \Big| 1 - \bE\big( e^{i\theta \Psi(X_1)/B_n}\big| \cF_{0}\big)\Big|^2 \to 0, \quad \theta \in \GR^1.\]
By Theorem \ref{PoCnew} it is enough to prove (\ref{eq:new2}), i.e.
\begin{align}
\Phi^h_n(\theta) &:= \sum_{j=1}^{n} \bE\big(e^{i\theta \Psi(X_{j})/B_n}\vert\cF_{j-1}\big) - 1  - i\theta B_n^{-1} \bE \big(\Psi(X_j) \Idn_{\{|\Psi(X_j)| \leq h B_n\}} \big| \cF_{j-1}\big)\label{eq:phien}\\
 &\inprob  \int \big( e^{i\theta x} - 1 - i \theta  x \Idn_{\{|x| \leq h\}}\big)\, \nu_{\alpha,c_+,c_-}(dx) =: \Phi^h(\theta).\nonumber
\end{align}
Let us notice that by (\ref{domattr1}) and  (\ref{domattr2}) we have
\[ \bE \Phi^h_n(\theta) \to \Phi^h(\theta).\]
Taking all these facts together we obtain the final condition to be verified:
\begin{equation}\label{eq:final}
\sum_{j=1}^{n} \chi^h_{n,\theta}(X_{j-1}) \inprob 0,
\end{equation}
where 
\begin{align*}
 \chi^h_{n,\theta}(x) &= \int \big(\exp\big(i\theta\Psi(y)/B_n\big) - 1 -i\theta \big(\Psi(y)/B_n\big) \Idn_{\{|\Psi(y)| \leq h B_n\}}\big)  P(x, dy) \nonumber \\
& - \Big(\int \big(\exp\big(i\theta\Psi(y)/B_n\big) - 1 -i\theta \big(\Psi(y)/B_n\big) \Idn_{\{|\Psi(y)| \leq h B_n\}} \big) \pi(dy)\Big).
\end{align*}
We apply the standard procedure based on the spectral gap property.
Let \[ \delta^h_{n,\theta} = \sum_{m=0}^{\infty} P^m(\chi^h_{n,\theta}).\]
$\delta^h_{n,\theta}$ is a well-defined element of $L^2_0(\pi)$, because $\chi^h_{n,\theta} \in L^2_0(\pi)$ and according to (\ref{eq:gap})
\begin{equation}\label{eq:delta}
\| \delta^h_{n,\theta} \|_2 \leq \sum_{m=0}^{\infty} \|P^m(\chi^h_{n,\theta})
\|_2 \leq \sum_{m=0}^{\infty} a^m \| \chi^h_{n,\theta}\|_2=\frac{\|\chi^h_{n,\theta}\|_2}{1-a}.
\end{equation}
Clearly,  $\delta^h_{n,\theta} = (I-P)^{-1}\chi^h_{n,\theta}$, hence $ \chi^h_{n,\theta} = \delta^h_{n,\theta}-P(\delta^h_{n,\theta})$.
Consequently, 
\[ \sum_{j=1}^{n} \chi^h_{n,\theta}(X_{j-1}) = \delta^h_{n,\theta}(X_0)-\delta^h_{n,\theta}(X_{n})+\sum_{j=1}^{n}\big(\delta^h_{n,\theta}(X_j)-P\delta^h_{n,\theta}(X_{j-1})\big).\]
The point here is that $\{\sum_{j=1}^{k}\big(\delta^h_{n,\theta}(X_j)-P\delta^h_{n,\theta}(X_{j-1})\big)\,;\, k=1,2,\ldots, n\}$ is a square integrable martingale. Therefore
\begin{align*} \bE\Big| \sum_{j=1}^{n} &\chi^h_{n,\theta}(X_{j-1})\Big|^2  \\&\leq 4\bE\vert\delta^h_{n,\theta}(X_0)\vert^2+4\bE\vert \delta^h_{n,\theta}(X_{n})\vert^2+4\bE\vert\sum_{j=1}^{n}\delta^h_{n,\theta}(X_j)-P\delta^h_{n,\theta}(X_{j-1})\vert^2\\
&=4\|\delta^h_{n,\theta}(X_0)\|_2^2+4\|\delta^h_{n,\theta}(X_{n})\|_2^2+4n\|\delta^h_{n,\theta}(X_1)-P\delta^h_{n,\theta}(X_0)\|_2^2\\
&\leq 8\|\delta^h_{n,\theta}(X_0)\|_2^2 +16n\|\delta^h_{n,\theta}(X_1)\|_2^2
\leq \frac{24}{(1-a)^2} n \bE\big|\chi^h_{n,\theta}(X_0)\big|^2 .
\end{align*}
So we have to prove that
\begin{equation}\label{eq:varpfi}
n \bE\big|\chi^h_{n,\theta}(X_0)\big|^2 \to 0.
\end{equation}
Let us notice that  
\[\|\chi^h_{n,\theta}(X_0)\|_2^2 = \bE\big|\chi^h_{n,\theta}(X_0)\big|^2 = \text{$\mathbb{V}$ar}\big(W^h_{n,\theta}\big) \leq  \bE\big|W^h_{n,\theta}\big|^2,\]
where 
\[ W^h_{n,\theta} =  1 + i\theta\bE \Big(\frac{\Psi(X_1)}{B_n} \Idn_{\{|\Psi(X_1)| \leq h B_n\}}\Big| \cF_{0}\Big) - \bE\big( e^{i\theta \Psi(X_1)/B_n}\big| \cF_{0}\big).\]
By inspection of (\ref{eq:basic}) we see that
\[ n\big|W^h_{n,\theta}\big|^2 \leq \frac{1}{2} \theta^4 I^h_{n,1} + 8 I^h_{n,2} \inprob 0.\]
As before, this convergence can be strengthened to the convergence in $L^1$ by applying the  $2$-U.I. condition. Indeed, both sequences  
\[\{Z'_n =\sqrt{n}\big(\frac{\Psi(X_1)^2}{B_n^2} \Idn_{\{|\Psi(X_1)| \leq h B_n\}}\big)\} \text{ and } \{Z^{\prime\prime}_n = \sqrt{n} \Idn_{\{|\Psi(X_1)| > h B_n\}}\}\]
 are bounded in $L^2$ by arguments essentially identical to those used in (\ref{eq:bound}) in the proof of $L^2$-boundedness of $\{Z_n = \sqrt{n}  \big(1 - e^{i\theta \Psi(X_1)/B_n}\big)\}$. Our $L^1$- claim follows from the observation that 
\[ n |W^h_{n,\theta}|^2 \leq \frac{1}{2}\theta^4 \Big(\bE\big(Z'_n\big|\cF_0\big)\Big)^2 + 8 \Big(\bE\big(Z^{\prime\prime}_n\big|\cF_0\big)\Big)^2.\]
We have thus completed the proof of Proposition \ref{PropMain}.

\begin{proposition}\label{PropMainPrim}
	In the framework of Proposition \ref{PropMain} replace the  $2$-U.I. condition with the hyperboundedness and the $L^2$-spectral gap property with strong mixing at geometric rate (or  geometric ergodicity). Then (\ref{eq:cent}) holds again.  
\end{proposition}

\noindent{\sc Proof.}\ The hyperboundedness of $P$ implies the  $2$-U.I.condition  and therefore most of the arguments used in the proof of Proposition \ref{PropMain} remain unchanged. The only place we used the $L^2$-spectral gap property was relation (\ref{eq:final}). We will show that the strong mixing in geometric rate and the hyperboundedness {\em taken together} imply more than required, namely
\begin{equation}\label{eq:todo}
 \bE\Big| \sum_{j=0}^{n-1} \chi^h_{n,\theta}(X_{j})\Big|^2 \to 0.
 \end{equation}
 By (\ref{eq:expo}) there is a number $0 \leq \eta < 1$ such that 
\[ \Big|\bE \Big(\chi_{n,\theta}^h(X_i)  \overline{\chi_{n,\theta}^h(X_j)}\Big)\Big|  \leq 2\pi \eta^{|i-j|} \big(2 + \theta h\big)^2, \quad i,j = 0,1,\ldots, n-1.\]
Set $m_n = [n^{\delta}]$, for some $\delta > 0$, and note that 
$n \eta^{m_n} \to 0$.
It follows that 
\begin{align*}
 \sum_{0 \leq i,j \leq n -1\atop |i-j| > m_n} &\Big|\bE \Big(\chi_{n,\theta}^h(X_i)  \overline{\chi_{n,\theta}^h(X_j)}\Big)\Big| \\
 &\leq 2\pi \big(2 + \theta h\big)^2\sum_{0 \leq i,j \leq n -1\atop |i-j| > m_n} \eta^{|i-j|} \leq \frac{4\pi \big(2 + \theta h\big)^2}{1 -\eta} \big(n-m_n) \eta^{m_n}.
 \end{align*}
So we have to consider the remaining covariances only.
\begin{align*}
\Big|\sum_{0 \leq i,j \leq n -1\atop |i-j| \leq  [n^{\delta}]} &\bE \Big(\chi_{n,\theta}^h(X_i)  \overline{\chi_{n,\theta}^h(X_j)}\Big)\Big| \leq  \sum_{i=0}^{n-1} \bE \Big|\chi_{n,\theta}^h(X_i)\Big|^2 \\
& + \sum_{i=0}^{n-2} \sum_{j=i+1}^{(i+[n^{\delta}])\wedge (n-1)} \Big|\bE \Big(\chi_{n,\theta}^h(X_i)  \overline{\chi_{n,\theta}^h(X_j)}\Big) + \bE \Big(\chi_{n,\theta}^h(X_j)  \overline{\chi_{n,\theta}^h(X_i)}\Big)\Big|\\
&\leq  n \bE \Big|\chi_{n,\theta}^h(X_0)\Big|^2 + \sum_{i=0}^{n-2} \sum_{j=i+1}^{(i+[n^{\delta}])\wedge (n-1)}
\bE \Big|\chi_{n,\theta}^h(X_i)\Big|^2  + \bE \Big|\chi_{n,\theta}^h(X_j)\Big|^2 \\
&\leq \big(n + 2(n-1)[n^{\delta}]\big) \bE \Big|\chi_{n,\theta}^h(X_0)\Big|^2 \leq 3 n^{1 + \delta}\ \bE \Big|\chi_{n,\theta}^h(X_0)\Big|^2,
\end{align*}
where the last inequality holds for $n$ large enough. It follows that the proof of Proposition \ref{PropMainPrim}
will be complete if we are able to prove that for some $\delta > 0$
\begin{equation}\label{eq:aim}
n^{1 + \delta}\ \bE \Big|\chi_{n,\theta}^h(X_0)\Big|^2 \to 0.
\end{equation}
Using the notation introduced in (\ref{eq:basic}) we have the following estimate.
\[n^{1+\delta} E \Big|\chi_{n,\theta}^h(X_0)\Big|^2 \leq 
n^{\delta} \Big( \frac{1}{2}\theta^4 \bE I^h_{n,1} + 8 \bE I_{n,2}^h \Big). \]
Notice that by Remark \ref{rem:rate} one can always find 
$\delta > 0$ such that 
$n^{\delta} \big( I^h_{n,1} + I^h_{n,2}\big) \inprob 0$.  
We shall strengthen this convergence using the hyperboundedness.
By this assumption there exists $q > 2$ such that the transition operator $P$ is a bounded linear map from $L^2(\pi)$ to $L^q(\pi)$: 
\[ \sup \{ \int \pi(dx) |Pf(x)|^{q}\,;\, \|f\|_2 \leq 1\} = \|P\|^{q}_{2\to q} < +\infty.\]
We have 
\begin{align*}
 \sup_n \Big\| \sqrt{n} &\frac{\Psi(X_1)^2}{B^2_n} \Idn_{\{|\Psi(X_1)| \leq h B_n\}}\Big\|^2_{2} = \\
 &=\sup_n n\ \bE
\Big(\frac{\Psi(X_1)^4}{B_n^4} \Idn_{\{|\Psi(X_1)| \leq h B_n\}}\Big) \\
&\leq h^2 \sup_n n\ \bE
\Big(\frac{\Psi(X_1)^2}{B_n^2} \Idn_{\{|\Psi(X_1)| \leq h B_n\}}\Big) \leq K_1 <+\infty.
\end{align*}
Therefore
\[\sup_n \bE \Big|\sqrt{n}\ \bE\Big(\frac{\Psi(X_1)^2}{B^2_n} \Idn_{\{|\Psi(X_1)| \leq h B_n\}}\Big| \cF_0\Big) \Big|^{q} \leq \|P\|^{q}_{2\to q} K_1^{q/2} < +\infty. \]  
Next we apply the H\"older inequality.
\begin{align*}
n^{\delta} \bE I_{n,1}^h &= n^{1 + \delta}\, \bE \Big(\bE\Big(\frac{\Psi(X_1)^2}{B^2_n} \Idn_{\{|\Psi(X_1)| \leq h B_n\}}\Big| \cF_0\Big)\Big)^2 \\
& = \bE \Bigg(\Big( \sqrt{n}\, \bE\Big(\frac{\Psi(X_1)^2}{B^2_n} \Idn_{\{|\Psi(X_1)| \leq h B_n\}}\Big| \cF_0\Big) \Big) \times \\
&\qquad\qquad\qquad \times \Big(n^{1/2 + \delta} \frac{\Psi(X_1)^2}{B^2_n} \Idn_{\{|\Psi(X_1)| \leq h B_n\}}\Big)\Bigg)\\
&\leq \Big(\bE \Big|\sqrt{n}\, \bE\Big(\frac{\Psi(X_1)^2}{B^2_n} \Idn_{\{|\Psi(X_1)| \leq h B_n\}}\Big| \cF_0\Big) \Big|^{q}\Big)^{1/q} \times \\
& \qquad\qquad\qquad \times \Big(\bE \Big|n^{1/2+\delta}\frac{\Psi(X_1)^2}{B^2_n} \Idn_{\{|\Psi(X_1)| \leq h B_n\}} \Big|^{\frac{q}{q - 1}} \Big)^{\frac{q-1}{q}} \\
&\leq \|P\|_{2\to q} K_1^{1/2} \Big(n^{\frac{(1/2+\delta) q}{q - 1}}\bE \Big|\frac{\Psi(X_1)}{B_n}\Big|^{\frac{2 q}{q - 1}} \Idn_{\{|\Psi(X_1)| \leq h B_n\}}  \Big)^{\frac{q-1}{q}}\\	
& \leq \|P\|_{2\to q} K_1^{1/2} h^{\frac{2}{q}}  \Big(n^{\frac{(1/2+\delta)q}{q - 1}} \bE \frac{\Psi(X_1)^2}{B^2_n} \Idn_{\{|\Psi(X_1)| \leq h B_n\}}  \Big)^{\frac{q-1}{q}} \to 0,
\end{align*}
if $ 0 <\delta < \frac{1}{2} - \frac{1}{q}$.

Similarly we handle  the other convergence $n^{\delta} \bE I_{n,2} \to 0$. First, we see that
\[
\sup_n \Big\| \sqrt{n} \Idn_{\{|\Psi(X_1)| > h B_n\}}\Big\|^2_{2}
=\sup_n n\ \bP\big(|\Psi(X_1)| > h B_n\big) \leq K_2 < +\infty.
\]
This implies that
\[\sup_n \bE \Big|\sqrt{n}\ \bP\big(|\Psi(X_1)| > h B_n\big| \cF_0\big) \Big|^{q} \leq \|P\|^{q}_{2\to q} K_2^{q/2} < +\infty. \]  
Finally we obtain by the H\"older inequality and for $ 0 <\delta < \frac{1}{2} - \frac{1}{q}$ that
\begin{align*}
n^{\delta} \bE &I_{n,2}^h = n^{1 + \delta}\, \bE \Big(\bP\big(|\Psi(X_1)| > h B_n\big| \cF_0\big)\Big)^2 \\
& = \bE \Big( \sqrt{n}\, \bP\big( |\Psi(X_1)| > h B_n\big| \cF_0\big) \Big) \Big(n^{1/2 + \delta}\, \Idn_{\{|\Psi(X_1)| > h B_n\}}\Big)\\
&\leq \Big(\bE \Big|\sqrt{n}\, \bP\big(|\Psi(X_1)| > h B_n\big| \cF_0\big) \Big|^{q}\Big)^{1/q}  \Big(\bE \Big|n^{1/2+\delta} \Idn_{\{|\Psi(X_1)| > h B_n\}} \Big|^{\frac{q}{q - 1}} \Big)^{\frac{q-1}{q}} \\
&\leq \|P\|_{2\to q} K_2^{1/2} \Big(n^{\frac{(1/2+\delta)q}{q - 1}}\bP \big(|\Psi(X_1)| > h B_n\big)  \Big)^{\frac{q-1}{q}} \to 0.
\end{align*}
The proof of Proposition \ref{PropMainPrim} is complete.

\subsection{Proof of Theorem \ref{TheMain}}
Given (\ref{eq:cent}), i.e. 
\[\frac{ \sum_{j=1}^{n} \Psi(X_j) - \bE \big(\Psi(X_j) \Idn_{\{|\Psi(X_j)| \leq h B_n\}} \big| \cF_{j-1}\big)}{B_n} 
\indist \text{$c_h$-Poiss}(\alpha, c_+, c_-),\]
we shall apply classic Theorem 4.2 from \cite{Bill68} in a way suitable for each case 
$\alpha \in (0,1)$, $\alpha = 1$ or $\alpha \in (1,2)$.

The reasoning is equally simple for $\alpha \neq 1$ and is based on the direct half of Karamata's theorem \cite[Theorem 1.5.11, p. 28]{BGT87}).

For $\alpha \in (0,1)$ we shall show that 
\begin{equation}\label{eq:lessone}
\lim_{h\to 0} \limsup_{n\to\infty} \bE\Big|\sum_{j=1}^{n}\bE \Big(\frac{\Psi(X_j)}{B_n} \Idn_{\{|\Psi(X_j)| \leq h B_n\}} \Big| \cF_{j-1}\Big)\Big| = 0, 
\end{equation}
and that
\[ \text{$c_h$-Poiss}(\alpha, c_+, c_-) \Rightarrow \mu_{\alpha},\ \text{ as $h\to 0$}.\]
The latter relation holds because $\int |x| \Idn_{\{ |x| \leq 1\}} \, \nu_{\alpha,c+,c_-}(dx) < +\infty$  for $\alpha \in (0,1)$. In order to prove (\ref{eq:lessone}) we proceed also in the standard way.
\begin{align*}
\bE\Big|\sum_{j=1}^{n}&\bE \Big(\frac{\Psi(X_j)}{B_n} \Idn_{\{|\Psi(X_j)| \leq h B_n\}} \Big| \cF_{j-1}\Big)\Big| \leq 
 n \bE \Big(\Big|\frac{\Psi(X_1)}{B_n}\Big| \Idn_{\{|\Psi(X_1)| \leq h B_n\}}\Big) \\
 &\leq \frac{n}{B_n} \int_0^{h B_n}\bP\big( |\Psi(X_1)| > t)\,dt \to_{n\to\infty} (1-\alpha)^{-1} h^{1-\alpha}(c_+ + c_-) \to_{h\to 0} 0.
 \end{align*}

For $\alpha \in (1,2)$ we shall show that 
\begin{equation}\label{eq:moretwo}
\lim_{h\to \infty} \limsup_{n\to\infty} \bE\Big|\sum_{j=1}^{n}\bE \Big(\frac{\Psi(X_j)}{B_n} \Idn_{\{|\Psi(X_j)| \geq h B_n\}} \Big| \cF_{j-1}\Big)\Big| = 0, 
\end{equation}
and that
\[ \text{$c_h$-Poiss}(\alpha, c_+, c_-) \Rightarrow \mu_{\alpha},\ \text{ as $h\to \infty$}.\]
Here again the latter relation holds due to the fact that  
\[\int |x| \Idn_{\{ |x| \geq 1\}} \, \nu_{\alpha,c+,c_-}(dx) < +\infty,\]  if  $\alpha \in (1,2)$. And we have 
\begin{align*}
\bE\Big|\sum_{j=1}^{n}&\bE \Big(\frac{\Psi(X_j)}{B_n} \Idn_{\{|\Psi(X_j)| \geq h B_n\}} \Big| \cF_{j-1}\Big)\Big| \leq 
n \bE \Big(\Big|\frac{\Psi(X_1)}{B_n}\Big| \Idn_{\{|\Psi(X_1)| \geq h B_n\}}\Big) \\
&\leq \frac{n}{B_n} \Big(\int_{h B_n}^{\infty}\bP\big( |\Psi(X_1)| > t)\,dt + h B_n \bP \big( |\Psi(X_1)| > h B_n\big)  \Big) \\ &\to_{n\to\infty} \frac{\alpha}{\alpha-1} h^{1-\alpha}(c_+ + c_-) \to_{h\to \infty} 0.
\end{align*}

The proof for $\alpha =1$ is somewhat different. Let us notice first that due to the symmetry of $\nu_{1,c,c}$ we have the equality
\[ \text{$c_{h}$-Poiss}(\alpha, c, c) =  \mu_{1}, \quad h \in \GR^1,\]
hence 
\[ \text{$c_{h}$-Poiss}(\alpha, c, c) \Rightarrow_{h\to 0} \mu_{1}.\]
Let $h > h' > 0$. By (\ref{eq:phien}) we have also
\[ \Phi^h_n(\theta) - \Phi^{h'}_{n}(\theta) = - i \theta \sum_{j=1}^{n} B_n^{-1} \bE \big(\Psi(X_j) \Idn_{\{h'B_n < |\Psi(X_j)| \leq h B_n\}} \big| \cF_{j-1}\big) \inprob 0, \ \ \theta \in \GR^1.\]
Therefore
\[ \lim_{h'\to 0} \limsup_{n\to\infty} \bP\Big( \Big|\sum_{j=1}^n \bE \Big(\frac{\Psi(X_j)}{B_n} \Idn_{\{h'B_n < |\Psi(X_j)| \leq h B_n\}} \Big| \cF_{j-1}\Big)\Big| > \varepsilon\Big) = 0, \ \ \varepsilon > 0. \]

\subsection{Proof of Theorem \ref{TheHyper}}
The starting point is the same as in the proof of Theorem
 \ref{TheMain}: by Proposition \ref{PropMainPrim} convergence (\ref{eq:cent}) holds and we have to reduce it to the desired form.
 The reduction for the cases $\alpha \in (0,1)$ and 
 $\alpha = 1$ is identical. So assume that $\alpha \in (1,2)$. Take $ 0 < r < \alpha -1$. We know that 
 \[ \int \pi(dx) |\Psi(x)|^{\alpha - r} < +\infty, \quad 0 < r <  \alpha - 1. \]
 We shall show that the hyperboundedness gives some $\beta > \alpha$ such that
 \[ \int \pi(dx) |P \Psi(x)|^{\beta} < +\infty.\]
 By the Riesz-Thorin interpolation theorem (see e.g. \cite[Theorem 1.1.1]{BeLo76}) applied to the transition  operator $P$ considered as a bounded linear map from $L^1(\pi)$ to $L^1(\pi)$ and from $L^2(\pi)$ to $L^{q}(\pi)$
 we have
 \[ \| P\|_{(\alpha - r) \to \beta} \leq \|P\|_{2\to  q}^{2(\alpha -r -1)/(\alpha -r)} < +\infty,\]
where 
\[ \beta = \frac{q(\alpha -r)}{2(q - 1) - (q - 2)(\alpha -r)} > \alpha. \]
Therefore 
\[ \bE \Big| \bE\big(\Psi(X_1)\big| \cF_0\big)\Big|^{\beta} < +\infty, \]
and we may apply Theorem \ref{TheWeak} to the conditional expectations  $\bE\big(\Psi(X_{j+1})\big| \cF_{j}\big)$, $j=0, 1,2,\ldots $ which are instantaneous functions of the strongly mixing at geometric rate Markov chain $\{X_j\}$.  

\begin{remark}\label{uzupelnienie}
In \cite{JKO09} a similar result on negligibility of sums of conditional expectations was obtained, but assuming the $L^2$-spectral gap property, which is stronger than the geometric ergodicity (and the strong mixing at geometric rate) and leads to a martingale decomposition (like in the proof of Theorem \ref{TheMain}).
	\end{remark}

\subsection{Proof of Theorem \ref{TheWeak}}
Let $\{X_n\}$ be strongly mixing at geometric rate and suppose that for some $\beta > 1$
\[ \int \pi(dx) |\Psi(x)|^{\beta} < +\infty,\quad \int \pi(dx) \Psi(x) = 0,\]
where $\pi$ is the stationary distribution.
Let $B_n = n^{1/\alpha}\tilde{\ell}(n)$, where $\tilde{\ell}(x)$ is a slowly varying function and 
$\alpha \in (0,\beta \wedge 2)$. Without loss of generality we may assume that $\beta <2$. 

Suppose we are able to prove that
\[ n \bE \Big| 1 - \bE\big( e^{i\theta \Psi(X_1)/B_n}\big| \cF_{0}\big)\Big|^2 = \bE\Big(\sum_{j=1}^{n} \Big| 1 - \bE\left(e^{i\theta \Psi(X_j)/B_n}\vert\cF_{j-1}\right) \Big|^2\Big) \to 0, \quad \theta \in \GR^1.\]
Then by Theorem \ref{PoCnew} it is enough to prove that
\[
\Phi_n(\theta) := \sum_{j=1}^{n} \bE\big(e^{i\theta \Psi(X_{j})/B_n}\vert\cF_{j-1}\big) - 1 \inprob 0, \theta \in \GR^1.
\]
Moreover, by the Marcinkiewicz-Zygmunt law of large numbers, we have $\bE \Phi_n(\theta) \to 0$ and so it is enough to prove that
\[
\bE\Big|\sum_{j=0}^{n-1} \big(\chi_{n,\theta}(X_{j}) - \bE \chi_{n,\theta}(X_{j}) \big)\Big|^2 \to 0,
\]
where 
\[ \chi_{n,\theta}(x) = \int \big(\exp\big(i\theta\Psi(y)/B_n\big) - 1 \big)  P(x, dy) - \Big(\int \big(\exp\big(i\theta\Psi(y)/B_n\big) - 1 \big) \pi(dy)\Big).\]
Similarly as in the proof of relation (\ref{eq:todo}), in presence of strong mixing at geometric rate it suffices to show that for some $\delta > 0$
\begin{equation}\label{eq:cruc}
 n^{1+\delta} \bE \big| \chi_{n,\theta}(X_{0}) \big|^2 \leq n^{1+\delta} \bE \Big| 1 - \bE\big( e^{i\theta \Psi(X_1)/B_n}\big| \cF_{0}\big)\Big|^2 \to 0.
 \end{equation}
 Notice that this convergence gives us also the crucial condition of Theorem \ref{PoCnew}. So proving (\ref{eq:cruc}) will complete the proof of Theorem \ref{TheWeak}.
 
 We have 
 \begin{align*}
 \Big| 1 - \bE\big( &e^{i\theta \Psi(X_1)/B_n}\big| \cF_{0}\big)\Big| \\
 &\leq |\theta| \bE \Big(\Big|\frac{\Psi(X_1)}{B_n}\Big| \Idn_{\{|\Psi(X_1)| \leq B_n\}} \Big| \cF_{0}\Big) + 2 \bP\big( |\Psi(X_1)| > B_n\big| \cF_0\big),
 \end{align*}
 and so, for $0 < \delta < \beta/\alpha - 1$ and $\beta <2$,
 \begin{align*}
n^{1 + \delta} \bE \Big| 1 - &\bE\big( e^{i\theta \Psi(X_1)/B_n}\big| \cF_{0}\big)\Big|^2 \\
&\leq 2|\theta|^2 n^{1 + \delta} \bE\Big|\bE \Big(\Big|\frac{\Psi(X_1)}{B_n}\Big| \Idn_{\{|\Psi(X_1)| \leq B_n\}} \Big| \cF_{0}\Big)\Big|^2\\
&\qquad\qquad + 4 n^{1 + \delta} \bE \Big(\bP\big( |\Psi(X_1)| > B_n\big| \cF_0\big)^2\Big) \\
&\leq 2|\theta|^2 n^{1 + \delta} \bE \Big|\frac{\Psi(X_1)}{B_n}\Big|^2 \Idn_{\{|\Psi(X_1)| \leq B_n\}} + 4 n^{1 + \delta} \bP\big(  |\Psi(X_1)| > B_n \big) \\
& \leq \big(2|\theta|^2 + 4\big)n^{1 + \delta} \bE \Big|\frac{\Psi(X_1)}{B_n}\Big|^{\beta}\\
& = \big(2|\theta|^2 + 4\big)n^{1 + \delta-\beta/\alpha} \big(\tilde{\ell}(n)\big)^{-\beta} \bE \big|\Psi(X_1)\big|^{\beta} \to 0.
 \end{align*}

\section*{Appendix: Complements on the Principle of Conditioning}\label{SecAppendix}

\setcounter{theorem}{0}
\renewcommand{\thetheorem}{A.\arabic{theorem}}

As mentioned in Introduction, the Principle of Conditioning (PoC) is a heuristic rule that allows producing limit theorems for dependent random variables given limit theorems for independent random variables. For example, applying the PoC  one obtains the following theorem on convergence to stable laws.  

\begin{theorem}\label{PoC}
Let $\{X_{n,j}\,;\,j\in \GN, n\in\GN\}$ be an array of random variables, which are row-wise adapted to
a sequence of filtrations $\{\{\mathcal{F}_{n,j}\,;\,j=0,1,\ldots\}\,;\,n\in\GN\}$. Let $h > 0$ and let $k_n\to\infty$ be a sequence of numbers.

The following conditions
\begin{align*}
\max_{1 \le j \le k_n} \bP \big( |X_{n,j}| > \varepsilon \big| \cF_{n,j-1}\big) &\inprob 0,\quad \varepsilon > 0;\\
\sum_{j=1}^{k_n} \bP\big( X_{n,j} > x| \cF_{n,j-1}\big) &\inprob c_+ x^{-\alpha},\quad x > 0;\\
\sum_{j=1}^{k_n} \bP\big( X_{n,j} < x\big| \cF_{n,j-1}\big) &\inprob c_- |x|^{-\alpha},\quad x < 0;\\
\sum_{j=1}^{k_n} \bE \big(X_{n,j}\Idnem_{\{|X_{n,j}| \leq h\}}\big| \cF_{n,j-1}\big) &\inprob a^h;\\
\sum_{j=1}^{k_n} \bV\text{ar} \big(X_{n,j}\Idnem_{\{|X_{n,j}| \leq h\}}\big| \cF_{n,j-1}\big) & \inprob \int_{\{|x| \leq h\}} x^2 \ \nu_{\alpha,c_+,c_-} (dx);
\end{align*}
imply that
\begin{equation}\label{eq:f2}
\sum_{j=1}^{k_n} X_{n,j} \indist \delta_{a^h}*\text{$c_h$-Poiss}(\alpha, c_+, c_-),
\end{equation}
where $\delta_{a^h}*\text{$c_h$-Poiss}(\alpha, c_+, c_-)$ is the stable distribution with the characteristic function (\ref{eq:f1}).
\end{theorem}
In other words the PoC says that if we replace in a limit theorem for row-wise independent summands:
\begin{itemize}
\item the expectations by conditional expectations with respect to the past,
\item the convergence of numbers by convergence in probability of random variables appearing in the conditions,
\end{itemize}
then  still the conclusion (in our case:  (\ref{eq:f2}))  will hold. In fact, one can also replace the summation to constants by summation to stopping times. 

We refer to \cite{Jak86} for exposition of results related to various versions of the PoC, beginning with the  Brown-Eagleson martingale CLT \cite{BrEa71}, through
multidimensional \cite{Klo77}, \cite{BKS82} and functional \cite{DuRe78}, \cite{JKM82} limit theorems, up to the PoC in infinite dimensional Hilbert \cite{Jak80}, \cite{Jak88} and Banach \cite{Ros82} spaces. The ideas standing behind the PoC motivated further research devoted to so called decoupling inequalities, described in detail in the well-known books
by Kwapie\'n and Woyczy\'nski \cite{KwWo92} and de la Pe\~{n}a and Gin\'e \cite{DPGi99}.  It might be interesting to realize that the tools developed to cope with the PoC find unexpected applications even today \cite{PeTa07}, \cite{JaRi16}. 

Behind the verbal form of the PoC there is a result on convergence of conditional characteristic functions (see \cite{Jak80}).

\begin{theorem}\label{PoCold} Let the system $\{X_{n,j}, \cF_{n,j}\}$ be as in Theorem \ref{PoC}. If for some $z \in \GC, z \neq 0,$ we have 
\[
\phi_n(\theta)=\prod_{j=1}^{k_n}\bE\left(e^{i\theta X_{n,j}}\vert\cF_{n,j-1}\right) \inprob z,\]
then also 
\[ \bE \exp (i\theta \sum_{j=1}^{k_n} X_{n,j}) \conver z.\]
In particular, if for some probability measure $\mu$ on $\GR^1$ we have
\begin{equation}\label{eq:chfunct}
 \phi_n(\theta) \inprob \hat{\mu}(\theta), \quad \theta \in \GR^1,
\end{equation}
then 
\[ \sum_{j=1}^{k_n} X_{n,j} \indist \mu.\]
\end{theorem}

Mimicking the case of independent random variables one can prove that conditions obtained by the PoC imply  (\ref{eq:chfunct}). But in many cases this is not the most efficient way of applying the PoC. It was observed in \cite{Jak12} that for highly structured models we can often check (\ref{eq:chfunct}) directly and that going this way we can keep integrability requirements at the minimal possible level.

We extend the results of \cite{Jak80}, \cite{JaKl80} and \cite{Jak12}
in the following theorem that provides  a convenient tool in many cases of interest.

\begin{theorem}\label{PoCnew}
Let $\{X_{n,j}\,;\,j\in \GN, n\in\GN\}$ be an array of random variables, which are row-wise adapted to
a sequence of filtrations $\{\{\mathcal{F}_{n,j}\,;\,j=0,1,\ldots\}\,;\,n\in\GN\}$.

Suppose that the following condition holds. 
\begin{equation}\label{eq:new1}
 \sum_{j=1}^{k_n} \big| 1 - \bE\left(e^{i\theta X_{n,j}}\vert\cF_{n,j-1}\right) \big|^2 \inprob 0, \quad \theta \in \GR^1.
\end{equation}
Let $A_n$ be arbitrary random variables and  $\Phi(\theta) \in \GC$ be a constant for each $\theta \in \GR^1$.
The following conditions are equivalent:
\begin{align}\label{eq:new2}
\Big(\sum_{j=1}^{k_n} \big(\bE\left(e^{i\theta X_{n,j}}\vert\cF_{n,j-1}\right) - 1\big)\Big) - i\theta A_n &\inprob \Phi(\theta).\\
\Big(\prod_{j=1}^{k_n} \bE\left(e^{i\theta X_{n,j}}\vert\cF_{n,j-1}\right)\Big) e^{-i\theta A_n} &\inprob e^{\Phi(\theta)}. \label{eq:new3}
\end{align}
In either case we have also
\begin{equation}\label{eq:new4} \bE \exp (i\theta \big(\sum_{j=1}^{k_n} X_{n,j} - A_n\big)) \conver e^{\Phi(\theta)}.
\end{equation}

In particular, if $e^{\Phi(\theta)} = \hat{\mu}(\theta)$, $\theta \in \GR^1$, for some probability measure $\mu$, then either of conditions (\ref{eq:new2}) or 
(\ref{eq:new3}) imply
\[ \sum_{j=1}^{k_n} X_{n,j} - A_n \indist \mu.\]
\end{theorem}

\noindent{\sc Proof.}\ Set 
\begin{align*}
\phi_n(\theta) &= \prod_{j=1}^{k_n} \bE\left(e^{i\theta X_{n,j}}\vert\cF_{n,j-1}\right); \\
\Phi_n(\theta) &= \sum_{j=1}^{k_n} \big(\bE\left(e^{i\theta X_{n,j}}\vert\cF_{n,j-1}\right) - 1\big).
\end{align*}
 If $z\in \GC$ satisfies $|z| \leq 1$, then $| z - e^{z-1}| \leq 5 |z - 1|^2$. Hence we have 
\begin{align*}
| \phi_n(\theta) e^{-i\theta A_n} - &\exp \big(\Phi_n(\theta) -i\theta A_n\big)| = | \phi_n(\theta) - \exp \big(\Phi_n(\theta)\big)| \\
&\leq \sum_{j=1}^{k_n} | \bE\left(e^{i\theta X_{n,j}}\vert\cF_{n,j-1}\right) - \exp\big( \bE\left(e^{i\theta X_{n,j}}\vert\cF_{n,j-1}\right) - 1\big)|\\
&\leq 5 \sum_{j=1}^{k_n} \big|\bE\left(e^{i\theta X_{n,j}}\vert\cF_{n,j-1}\right) - 1\big|^2 \inprob 0, \quad\text{ by (\ref{eq:new1})}.
\end{align*}
We have thus established the equivalence of (\ref{eq:new2}) and
(\ref{eq:new3}). To prove that (\ref{eq:new3}) implies (\ref{eq:new4}) we need a suitable version of Lemma 2 in \cite{JaSl86}.

\begin{lemma}
For every $\varepsilon >0$
\begin{equation}
\begin{split}
 \big| \bE \exp (i\theta &\big(\sum_{j=1}^{k_n} X_{n,j} - A_n\big)) - \bE \big(\phi_n(\theta) e^{-i\theta A_n}\big)\big| \\
 &\leq 2(1+\frac{1}{\varepsilon}) \bP\big(|\phi_n(\theta)| < \varepsilon\big) + \frac{1}{\varepsilon} \bE \big| \phi_n(\theta) e^{-i\theta A_n} - \bE \big(\phi_n(\theta) e^{-i\theta A_n}\big)\big|.
\end{split} 
\end{equation}
\end{lemma}

\noindent{\sc Proof.}\ We follow the idea of the proof of Theorem A in \cite{Jak80}. Define
\[ \phi_{n,k}(\theta) = \prod_{j=1}^{k} \bE\left(e^{i\theta X_{n,j}}\vert\cF_{n,j-1}\right).\]
Fix $\theta \in \GR^1$ and $\varepsilon > 0$ and consider random variables
\[ X_{n,k}^* = X_{n,k} \Idn_{\{|\phi_{n,k}(\theta)| \geq \varepsilon\}}.\]
Then we have both
\[ \big| \bE \exp (i\theta \big(\sum_{j=1}^{k_n} X_{n,j} - A_n\big)) -
\bE \exp (i\theta \big(\sum_{j=1}^{k_n} X^*_{n,j} - A_n\big)) \big| \leq 2\bP\big( |\phi_n(\theta)| < \varepsilon\big),\]
and, if we set $\phi_n^*(\theta) = \prod_{j=1}^{k_n} \bE\left(e^{i\theta X^*_{n,j}}\vert\cF_{n,j-1}\right)$,
\[ 
 \bE \big| e^{-i\theta A_n} \phi_n(\theta) - e^{-i\theta A_n} \phi^*_n(\theta) \big|  \leq 2\bP\big( |\phi_n(\theta)| < \varepsilon\big).
\]
The advantage of random variables $\{X^*_{n,j}\}$ consists in the fact that
\[ |\phi_n^*(\theta)| = \big| \prod_{j=1}^{k_n} \bE\left(e^{i\theta X^*_{n,j}}\vert\cF_{n,j-1}\right) \big| \geq \varepsilon,\]
and so, by the backward induction (or the martingale property)
\[ \bE \frac{\exp (i\theta \big(\sum_{j=1}^{k_n} X^*_{n,j} - A_n\big))}{e^{-i\theta A_n}\phi_n^*(\theta)} = \bE \frac{\exp (i\theta \big(\sum_{j=1}^{k_n} X^*_{n,j}\big))}{\prod_{j=1}^{k_n} \bE\left(e^{i\theta X^*_{n,j}}\vert\cF_{n,j-1}\right)}  =1.\]
Therefore, 
\begin{align*}
\big|\bE \exp (i\theta &\big(\sum_{j=1}^{k_n} X^*_{n,j} - A_n\big)) - 
 \bE \big(\phi_n(\theta) e^{-i\theta A_n}\big)\big| \\
&= \Big|\bE \frac{\exp (i\theta \big(\sum_{j=1}^{k_n} X^*_{n,j}\big))}{\phi_n^*(\theta)} \phi_n^*(\theta) e^{-i\theta A_n} \\
 &\qquad\qquad - \bE \big(\phi_n(\theta) e^{-i\theta A_n}\big)\bE \frac{\exp (i\theta \big(\sum_{j=1}^{k_n} X^*_{n,j}\big))}{\phi_n^*(\theta)} \Big|\\
&\leq \frac{1}{\varepsilon} \bE |  \phi_n^*(\theta) e^{-i\theta A_n} - \bE  \phi_n(\theta) e^{-i\theta A_n}| \\
&\leq \frac{1}{\varepsilon} \Big( 2\bP\big( |\phi_n(\theta)| < \varepsilon\big) + \bE \big| \phi_n(\theta) e^{-i\theta A_n} - \bE \big(\phi_n(\theta) e^{-i\theta A_n}\big)\big| \Big).
\end{align*}

\noindent{\sc Proof of Theorem \ref{PoCnew} (continued).}
Now assume that (\ref{eq:new3}) holds. Let $\varepsilon = 1/2|e^{\Phi(\theta)}|$. Then $\bP\big( |\phi_n(\theta)| < \varepsilon\big) = \bP\big( |\phi_n(\theta)e^{-i\theta A_n}| < \varepsilon\big)
\to 0$ and by the dominated convergence $\bE \big| \phi_n(\theta) e^{-i\theta A_n} - \bE \big(\phi_n(\theta) e^{-i\theta A_n}\big)\big| \to 0$.

\section*{Acknowledgments}
The authors thank Zbigniew Szewczak for pointing us the article by Wu \cite{Wu00}.  
A. Jakubowski gratefully acknowledges the support of the University of Rouen Normandie.


\section*{References}

\begin{thebibliography}{10}
	
	\bibitem{BJMW11}
	K.~Bartkiewicz, A.~Jakubowski, T.~Mikosch, and O.~Wintenberger.
	\newblock Stable limits for sums of dependent infinite variance random
	variables.
	\newblock {\em Probab. Theory Relat. Fields}, 150:337--372, 2011.
	
	\bibitem{BeLo76}
	J.~Bergh and J.~L\"ofstr\"om.
	\newblock {\em Interpolation Spaces. \protect{A}n introduction}, volume 223 of
	{\em Grundlehren Math. Wiss.}
	\newblock Springer, Heidelberg, 1976.
	
	\bibitem{BKS82}
	M.~Be\'ska, A.~K\l opotowski, and L.~S\l omi\'nski.
	\newblock Limit theorems for random sums of dependent d-dimensional random
	vectors.
	\newblock {\em Z. Wahrscheinlichkeitstheor. Verw. Gebiete}, 61:43--57, 1982.
	
	\bibitem{Bill68}
	P.~Billingsley.
	\newblock {\em Convergence of Probability Measures}.
	\newblock Wiley, New York, 1968.
	
	\bibitem{BGT87}
	N.H. Bingham, C.M. Goldie, and J.L. Teugels.
	\newblock {\em Regular Variation}, volume~27 of {\em Encyclopedia Math. Appl.}
	\newblock Cambridge Univ. Press, Cambridge, 1987.
	
	\bibitem{Bra83}
	R.C. Bradley.
	\newblock Information regularity and the central limit question.
	\newblock {\em Rocky Mountain J. Math.}, 13:77--97, 1983.
	
	\bibitem{Bra07I}
	R.C. Bradley.
	\newblock {\em Introduction to Strong Mixing Conditions, Vol. I}.
	\newblock Kendrick Press, Heber City, 2007.
	
	\bibitem{Bra07II}
	R.C. Bradley.
	\newblock {\em Introduction to Strong Mixing Conditions, Vol. II}.
	\newblock Kendrick Press, Heber City, 2007.
	
	\bibitem{BrEa71}
	B.~M. Brown and G.~K. Eagleson.
	\newblock Martingale convergence to infinitely divisible laws with finite
	variances.
	\newblock {\em Trans. Amer. Math. Soc.}, 162:449--453, 1971.
	
	\bibitem{BDM16}
	D.~Buraczewski, E.~Damek, and T.~Mikosch.
	\newblock {\em Stochastic Models with Power-Law Tails. The Equation $X = AX +
		B$}.
	\newblock Springer, 2016.
	
	\bibitem{CaM-A14}
	P.~Cattiaux and M.~Manou-Abi.
	\newblock Limit theorems for some functionals with heavy tails of a discrete
	time markov chain.
	\newblock {\em ESAIM Probab. Stat.}, 18:468--486, 2014.
	
	\bibitem{Cog72}
	R.~Cogburn.
	\newblock The central limit theorem for \protect{Markov} processes.
	\newblock In L.M.~Le Cam, J.~Neyman, and E.L. Scott, editors, {\em Proceedings
		of the Sixth Berkeley Symposium on Mathematical Statistics and Probability,
		Volume 2: Probability Theory}, pages 485--512. University of California
	Press, 1972.
	
	\bibitem{Dav83}
	R.A. Davis.
	\newblock Stable limits for partial sums of dependent random variables.
	\newblock {\em Ann. Probab.}, 11:262--269, 1983.
	
	\bibitem{DaMi98}
	R.A. Davis and T.~Mikosch.
	\newblock Limit theory for the sample acf of stationary process with heavy
	tails with applications to ARCH.
	\newblock {\em Ann. Probab.}, 26:2049--2080, 1998.
	
	\bibitem{Davy73}
	Yu.~A. Davydov.
	\newblock Mixing conditions for Markov chains.
	\newblock {\em Teor. Veroyatnost. i Primienen.}, 18:321--338, 1973.
	
	\bibitem{DPGi99}
	V.~de~la Pe\~{n}a and E.~Gin\'e.
	\newblock {\em Decoupling: from dependence to independence}.
	\newblock Springer, New York, 1999.
	
	\bibitem{DeJa89}
	M.~Denker and A.~Jakubowski.
	\newblock Stable limit distributions for strongly mixing sequences.
	\newblock {\em Statist. Probab. Lett.}, 8:477--483, 1989.
	
	\bibitem{Doob53}
	J.L. Doob.
	\newblock {\em Stochastic Processes}.
	\newblock Wiley, New York, 1953.
	
	\bibitem{DuRe78}
	R.~Durret and S.~Resnick.
	\newblock Functional limit theorems for dependent random variables.
	\newblock {\em Ann. Probab.}, 6:829--846, 1978.
	
	\bibitem{EKM97}
	P.~Embrechts, C.~Kl\"{u}ppelberg, and T.~Mikosch.
	\newblock {\em Modelling Extremal Events for Insurance and Finance}.
	\newblock Springer, Berlin, 1997.
	
	\bibitem{Fell70}
	W.~Feller.
	\newblock {\em An Introduction to Probability Theory and Its Applications.
		Volume II. Second Edition}.
	\newblock Wiley, New York, 1970.
	
	\bibitem{GoLi78}
	M.I. Gordin and B.A. Lifshitz.
	\newblock A central limit theorem for Markov processes.
	\newblock {\em Soviet Math. Doklady}, 19:392--394, 1978.
	
	\bibitem{Hagg05}
	O.~H\"aggstr\"om.
	\newblock On the central limit theorem for geometrically ergodic Markov chains.
	\newblock {\em Probab. Theory Relat. Fields}, 132:74--82, 2005.
	
	\bibitem{Hagg06}
	O.~H\"aggstr\"om.
	\newblock Acknowledgement of priority concerning ``On the central limit theorem
	for geometrically ergodic Markov chains".
	\newblock {\em Probab. Theory Relat. Fields}, 135:470, 2006.
	
	\bibitem{JKM82}
	J.~Jacod, A.~K\l opotowski, and J.~M\'emin.
	\newblock Th\'eor\`eme de la limite centrale et convergence fonctionelle vers
	un processus \`a accroissements ind\'ependants: la m\'ethode des martingales.
	\newblock {\em Ann. Inst. H. Poincar\'e, Sect. B}, 18:1--45, 1982.
	
	\bibitem{JaSh03}
	J.~Jacod and A.~Shiryayev.
	\newblock {\em Limit theorems for stochastic processes. \protect{S}econd
		edition}, volume 288 of {\em Grundlehren Math. Wiss.}
	\newblock Springer, Heidelberg, 2003.
	
	\bibitem{Jak80}
	A.~Jakubowski.
	\newblock On limit theorems for sums of dependent Hilbert space valued random
	variables.
	\newblock {\em Lecture Notes in Statist.}, 2:178--187, 1980.
	
	\bibitem{Jak86}
	A.~Jakubowski.
	\newblock Principle of conditioning in limit theorems for sums of random
	variables.
	\newblock {\em Ann. Probab.}, 14:902--915, 1986.
	
	\bibitem{Jak88}
	A.~Jakubowski.
	\newblock Tightness criteria for random measures with application to the
	principle of conditioning in Hilbert spaces.
	\newblock {\em Probab. Math. Statist.}, 9:95--114, 1988.
	
	\bibitem{Jak12}
	A.~Jakubowski.
	\newblock Principle of conditioning revisited.
	\newblock {\em Demonstratio Math.}, XLV:325--336, 2012.
	
	
	
	\bibitem{JaKl80}
	A.~Jakubowski and A.~K\l opotowski.
	\newblock Quelques remarques sur les d\'emonstrations des th\'eor\`emes limite
	pour des vecteurs $d$-dimensionells al\'eatoires non ind\'ependants.
	\newblock {\em Publ. S\'eminaire de Probabilit\'es de Rennes}, exp. no 3:1--16,
	1980.
	
	\bibitem{JaRi16}
	A.~Jakubowski and M.~Riedle.
	\newblock Stochastic integration with respect to cylindrical L\'evy processes.
	\newblock {\em Ann. Probab.}, 45:4273--4306, 2018.
	
	\bibitem{JaSl86}
	A.~Jakubowski and L.~S\l omi\'nski.
	\newblock Extended convergence to continuous in probability processes with
	independent increments.
	\newblock {\em Probab. Theory Related Fields}, 72:55--82, 1986.
	
	\bibitem{JKO09}
	M.~Jara, T.~Komorowski, and S.~Olla.
	\newblock Limit theorems for additive functionals of a \protect{M}arkov chain.
	\newblock {\em Ann. Appl. Probab.}, 19:2270--2300, 2009.
	
	\bibitem{Klo77}
	A.~K\l opotowski.
	\newblock Limit theorems for sums of dependent random vectors in
	$\mathbb{R}^d$.
	\newblock {\em Dissertationes Math.}, 151:1--62, 1977.
	
	\bibitem{KoMe12}
	I.~Kontoyiannis and S.~Meyn.
	\newblock Geometric ergodicity and the spectral gap of non-reversible
	\protect{M}arkov chains.
	\newblock {\em Probab. Theory Relat. Fields}, 154:327--339, 2012.
	
	\bibitem{Kriz10}
	D.~Krizmanic.
	\newblock {\em Functional limit theorems for weakly dependent regularly varying
		time series}.
	\newblock PhD thesis, University of Zagreb, 2010.
	
	\bibitem{KwWo92}
	S.~Kwapie\'n and W.~A. Woyczy\'nski.
	\newblock {\em Random Series and Stochastic Integrals: Single and Multiple}.
	\newblock Birkh\"auser, Basel, 1992.
	
	\bibitem{MaWo00}
	M.~Maxwell and M.~Woodroofe.
	\newblock Central limit theorems for additive functionals of Markov chains.
	\newblock {\em Ann. Probab.}, 28:713--724, 2000.
	
	\bibitem{MeTw09}
	S.~Meyn and R.L. Tweedie.
	\newblock {\em Markov Chains and Stochastic Stability. Second Edition}.
	\newblock Cambridge, Cambridge, 2009.
	

	
	\bibitem{PeTa07}
	G.~Peccati and M.S. Taqqu.
	\newblock Stable convergence of generalized $l^2$ stochastic integrals and the
	principle of conditioning.
	\newblock {\em Electron. J. Probab.}, 12:447--480, 2007.
	
	\bibitem{RoRo97}
	G.O. Roberts and J.S. Rosenthal.
	\newblock Geometric ergodicity and hybrid Markov chains.
	\newblock {\em Electron. Commun. Probab.}, 2:13--25, 1997.
	
	\bibitem{Ros71}
	M.~Rosenblatt.
	\newblock {\em Markov Processes. Structure and Asymptotic Behavior}.
	\newblock Springer Verlag, New York, 1971.
	
	\bibitem{Ros82}
	J.~Rosi\'nski.
	\newblock Central limit theorems for dependent random vectors in
	\protect{B}anach spaces, in: \protect{J-A. Chao, W.A. Woyczy\'nski, Eds.},
	\protect{M}artingale theory in harmonic analysis and \protect{B}anach spaces.
	\newblock {\em Lecture Notes in Math.}, 939:157--180, 1982.
	
	\bibitem{SaTa94}
	G.~Samorodnitsky and M.S. Taqqu.
	\newblock {\em Stable Non-Gaussian Random Processes: Stochastic Models with
		Infinite Variance}.
	\newblock Chapman \& Hall/CRC, Boca Raton, 1994.
	
	\bibitem{Wu00}
	L.~Wu.
	\newblock Uniformly integrable operators and large deviations for Markov
	processes.
	\newblock {\em J. Funct. Anal.}, 172:301--376, 2000.
	
\end{thebibliography}

\end{document}